\documentclass[11pt]{amsart}
\title{Realization of Cohomology Classes in Grassmannians}
\author{Izzet Coskun and Julius Ross}
\date{\today}
\usepackage{amsmath,amssymb,amsthm, latexsym, amsopn, array}
\usepackage[colorlinks,linkcolor=blue]{hyperref}
\usepackage{xcolor}

\theoremstyle{plain}
\newtheorem{lemma}{Lemma}[section]
\newtheorem*{theorem*}{Theorem}
\newtheorem*{lemma*}{Lemma}
\newtheorem*{proposition*}{Proposition}
\newtheorem*{conjecture*}{Conjecture}
\newtheorem*{corollary*}{Corollary}
\newtheorem*{problem*}{Problem}
\newtheorem{theorem}[lemma]{Theorem}

\newtheorem{corollary}[lemma]{Corollary}
\newtheorem{proposition}[lemma]{Proposition}

\newtheorem{problem}[lemma]{Problem}
\newtheorem{question}[lemma]{Question}

\theoremstyle{definition}
\newtheorem{definition}[lemma]{Definition}

\newtheorem{remark}[lemma]{Remark}

\begin{document}

\newcommand{\Hom}{\operatorname{Hom}}
\newcommand{\GG}{\mathbb{G}}
\newcommand{\PP}{\mathbb{P}}
\newcommand{\ZZ}{\mathbb{Z}}
\newcommand{\QQ}{\mathbb{Q}}
\newcommand{\CC}{\mathbb{C}}
\newcommand{\OO}{\mathcal{O}}

\newcommand{\Span}{\operatorname{Span}}
\begin{abstract}
The study of irreducible subvarieties has recently seen a surge of interest due to connections with convex geometry. In this paper, we study cohomology classes of Grassmannians that are realizable by irreducible subvarieties.  We completely classify the cohomology classes that can be realized by irreducible subvarieties in dimensions 2 and 3 and codimensions 2 and 3. We also classify all the classes that can be realizable by an irreducible subvariety for the Grassmannians $G(2,n)$ for $n \leq 6$ and $G(3,6)$. We classify the cohomology classes in all dimensions that are, up to positive multiple, realizable by an irreducible subvariety for the Grassmannians $G(2,n)$.
\end{abstract}

\maketitle
\section{Introduction}

Let $G(k,n)$ denote the Grassmannian parameterizing $k$-dimensional subspaces of an $n$-dimensional vector space.  Given a subvariety $X \subset G(k,n)$,  we define the cohomology class $[X]$ of $X$ to be the Poincar\'e dual of its fundamental class. In this paper, we address the following problem.

\begin{problem}\label{problem-main}
    When is a cohomology class $\nu \in H^{2m}(G(k,n), \ZZ)$  the class of an irreducible subvariety of $G(k,n)$?
\end{problem}

Following \cite{hhmww}, we say that the class $\nu$ is {\em realizable over $\ZZ$} if $\nu$ is the class of an irreducible subvariety of $G(k,n)$. Similarly, the class $\nu$ is {\em realizable over $\QQ$} if some positive rational multiple of $\nu$ is the class of an irreducible subvariety. 

The cohomology of $G(k,n)$ has an additive basis given by Schubert classes $\sigma_{\lambda}$, indexed by partitions $\lambda$ with at most $k$ parts and whose parts are at most $n-k$. More explicitly, Problem \ref{problem-main} asks for which integral linear combinations $\nu = \sum a_{\lambda} \sigma_{\lambda}$ of Schubert classes does there exists an irreducible subvariety $X$ of $G(k,n)$ such that $[X]=\nu$? By Kleiman transversality, any effective class is a nonnegative linear combination of Schubert cycles (see, for example, \cite[Proposition 2.20]{coskun:IMPANGA}). Hence, a necessary condition is that $a_{\lambda} \geq 0$. 

It is well-known that in dimension 1 and in codimension 1, every effective class can be represented by an irreducible subvariety of $G(k,n)$. Similarly, in dimension 2 and in codimension 2, with the exception of multiples of Schubert classes, every effective class can be represented by irreducible subvarieties. 

\begin{theorem}[Theorem \ref{thm-dimension2} and Theorem \ref{thm-codimension2}] 
  Let $\nu = \sum a_{\lambda} \sigma_{\lambda}$ be a dimension 2 or codimension 2 cohomology class in $G(k,n)$. If $a_{\lambda} >0$ for every $\lambda$, then $\nu$ is realizable over $\ZZ$.  
\end{theorem}
This result on dimension 2 classes was recently also obtained in \cite[Theorem 7.7]{hhmww}.

Starting in dimension and codimension 3, there are additional restrictions on the classes that are realizable coming from the Hodge-Riemann relations. In dimension 3 and codimension 3, these are essentially the only restrictions. 

\begin{theorem}[Theorem \ref{thm-dimension3} and Theorem \ref{thm-codimension3}]
Let $3 \leq k \leq n-3$ be integers. 
    Let $$\nu = a \sigma_{3} + b \sigma_{2,1} + c \sigma_{1,1,1}$$ or $$\nu = a \sigma_{(n-k)^{k-1}, n-k-3} + b \sigma_{(n-k)^{k-2}, n-k-2, n-k-1}+ c \sigma_{(n-k)^{k-3}, (n-k-1)^3}$$ be a codimension 3 or, respectively, dimension 3 cohomology class in $G(k,n)$ with $a,b,c>0$. Then $\nu$ is realizable over $\ZZ$ if and only if $b^2 \geq ac$.
\end{theorem}

Similarly, we classify all the cohomology classes that are realizable over $\ZZ$ for the Grassmannians $G(2,n)$ for $n \leq 6$ (see Theorem \ref{thm-dim4G(2,n)}) and for the Grassmannian $G(3,6)$ (see Theorem \ref{thm-5dimG36}).

In general, classifying which cohomology classes are realizable over $\ZZ$ is a difficult problem. There exist classes for which some integer multiples are realizable over $\ZZ$, while others are not. For example, although Schubert classes themselves are realizable over 
$\ZZ$, their multiples may fail to be. The problem becomes more tractable over $\QQ$. We classify all the classes which are realizable over $\QQ$ in Grassmannians $G(2,n)$.

\begin{theorem}[Theorem \ref{thm-g2ningeneral}]
    The cohomology class $\sum a_i \sigma_{n-2-i, i}$ is realizable over $\QQ$ in $G(2,n)$ if and only if the coefficients $a_i$ form a log concave sequence of nonnegative  rational numbers with no internal zeros.
\end{theorem}

We do not know a complete classification of realizable classes in higher dimensions and codimensions in general. We do, however, give constructions of irreducible subvarieties via an inductive process. We also give a general method of finding obstructions to realizability coming from the Hodge-Riemann relations. Moreover, we prove that the realizable classes stabilize as $n$ increases.

\begin{theorem}[Corollaries \ref{cor-equalitiesIJ} and \ref{cor-equalitydim}]\label{thm:stabilization:intro} 
    Let $k,n,r$ be positive integers such that $n \geq k+r$. Then a cohomology class 
    $\nu = \sum a_{\lambda} \sigma_{\lambda}$ of dimension or codimension $r$ is realizable over $\QQ$ in $G(k,n)$ if and only if it is realizable over $\QQ$ in $G(k, k+r)$. Moreover, if $r \leq k \leq n-r$, then $\nu$ is realizable over $\QQ$ in $G(k,n)$ if and only if it is realizable over $\QQ$ in $G(r, 2r)$
\end{theorem}

In fact, we will see that to understand codimension $r$ classes that are realizable over   $\ZZ$ in a Grassmannian $G(k,n)$ with $k \leq r \leq n-r$, it suffices to understand those that are realizable in $G(r, 2r)$. We note, however, that over $\ZZ$ stabilization occurs only starting with $G(r+1, 2r+2)$ (see \S \ref{sec:reverseinclusions}).

\subsection*{Comparsion with Other Works}

Problem \ref{problem-main} has been studied extensively for multiples of Schubert classes in homogeneous varieties (\cite{bryant:rigidity, coskun:rigid, coskun:rigidorth, coskun:IMPANGA, coskunrobles,  hong:rigidity, hong:singular, HongMok, liu1, Liu, Liu2, RoblesThe, walters:thesis}). There is a complete classification of multiples of Schubert classes that are realizable over $\ZZ$ in compact complex Hermitian symmetric spaces \cite{coskunrobles, RoblesThe}. We will recall the precise classification in the case of $G(k,n)$ in \S \ref{sec-prelim} since these will play a crucial role in the proof of our classification. 

June Huh has made significant progress in the related problem of understanding realizable classes in products of projective spaces \cite{huh, huh2}. Recently, Huang, Huh, Micha\l{}ek, Wang and Wang classified surfaces that are realizable over $\QQ$ in products of projective spaces and surfaces that are realizable over $\ZZ$ in Grassmannians \cite{hhmww}.

\subsection*{Questions and Further Directions}

Our first question concerns a duality between dimension and codimension.

\begin{question}
Does the duality given by taking the complement of a partition preserve being realizable over $\QQ$?  Said another way, if $\lambda^c$ denotes the complement of a partition $r$, is it the case that a dimension $r$ cohomology class $\sum_{\lambda} a_{\lambda} \sigma_{\lambda}$  is realizable over $\QQ$ if and only if the codimension $r$ cohomology class
$\sum_{\lambda} a_{\lambda} \sigma_{\lambda^c}$ is realizable over $\QQ$?
\end{question}

Such duality holds in all the cases we have computed.  From Theorem \ref{thm-codimension3} we see this duality holding between dimension 3 and codimension 3, but we are not aware of any way to directly get the corresponding representative of a dimension 3 class from a representative from the dual codimension 3 class.

Our second question comes from the obstructions we obtain to realizability, all of which come from the Hodge-Riemann relations.  Given a codimension $r$ class $\nu=\sum_{\lambda} a_{\lambda} \sigma_{\lambda}$  we know from Theorem \ref{thm:stabilization:intro} that  $\eta$ is realizable in $G(k,k+r)$ if and only if it is realizable in $G(k,n)$ for all $n\ge k+r$.   For such an $n$ suppose we can find a partition $k_1+\cdots+k_p=k$ and further positive integers $n_1,\ldots,n_p$ such that the natural rational map
\begin{equation}\phi: G(k_1,n_1)\times G(k_2,n_2)\times \cdots \times G(k_p,n_p)\dashrightarrow G(k,k+n)\label{eq:productsofgrassmannians}\end{equation}
obtained by taking the span of the elements on the left hand side is birational.     Then pulling back a translate of an irreducible representative of $\nu$ gives an irreducible representative $Y$ of $\phi^*\nu$.    Letting $H_1,\ldots,H_p$ denote the hyperplane class of the factors on the left hand side, the Hodge-Riemann relations imply that if $\alpha_1+ \cdots + \alpha_p = \dim(Y)-2$ then the matrix $L=(L_{uv})$ given by
\begin{equation*}L_{uv}= \int_Y H_1^{\alpha_1}\cdots H_p^{\alpha_p} H_u H_v \text{ for }u,v=1,\ldots,p\label{eq:thematrixL}\end{equation*}
is weakly Lorentzian (i.e.,\ it is the limit of matrices that have signature $(+,-,\cdots,-)$).

\begin{question}
Does this process give all obstructions to a cohomology class in the Grassmannian being realizable over $\QQ$?  Specifically, is it the case that $\nu$ is realizable over $\QQ$ if and only if the matrix $L$ is weakly Lorentzian for all $k\ge r$ all $n\ge k+r$ and all such $k_i,n_i,\alpha_i$? If not, what other restrictions are there?
\end{question}

One may ask a similar question for dimension $r$ classes.  By the work in this paper we know the answer is yes for dimension 3 classes and for codimension 3 class, and is trivially true for dimension 2 and codimension 2 classes.  If the above has a positive answer it would be interesting to know if the set of such obstructions can be explicitly described by a finite number of conditions on the coefficients of $\nu$.

\subsection*{Acknowledgements} IC is partially supported by the National Science Foundation grant DMS 2200684. JR is partially supported by the National Science Foundation grant DMS 1749447.  Both authors are also supported by a Simons Foundation Award.   We would like to thank Olivier Debarre for sparking our interest in this problem and generously sharing his ideas.

\section{Preliminaries}\label{sec-prelim}

Let $V_n$ be an $n$-dimensional complex vector space. We let $G(k, V_n)$ denote the Grassmannian parameterizing $k$-dimensional subspaces of $V_n$. When the subspace $V_n$ is immaterial, we will denote this Grassmannian by $G(k,n)$.  Let $\GG(k-1, \PP V_n)$ denote the Grassmannian parameterizing $(k-1)$-dimensional projective linear spaces in $\PP V_n$. We will often denote this Grassmannian by $\GG(k-1, n-1)$. 
The dimension of $G(k,n)$ is $k(n-k)$.

Let $\lambda=(\lambda_1,\ldots,\lambda_k)$ be a partition with $k$ parts satisfying $$n-k \geq \lambda_1 \geq \cdots\ge\lambda_k \geq 0.$$ Let $$F_{\bullet} : F_1 \subset \cdots \subset F_n = V_n$$ be a complete flag in $V_n$, where $\dim (F_i) = i$. Then the Schubert variety $\Sigma_{\lambda}(F_{\bullet}) \subset G(k,n)$ is defined by
$$\Sigma_{\lambda}(F_{\bullet}) := \{ W \in G(k,n) | \dim(W \cap F_{n-k+i - \lambda_i}) \geq i \ \mbox{for} \ 1 \leq i \leq k \}.$$ 
The codimension of $\Sigma_{\lambda}(F_{\bullet})$ is $|\lambda| = \sum_{i=1}^k \lambda_i$, and the dimension of $\Sigma_{\lambda}(F_{\bullet})$ is $k(n-k) - |\lambda|.$

We denote the Poincar\'{e} dual of $\Sigma_{\lambda}(F_{\bullet})$ by $\sigma_{\lambda} \in H^{2|\lambda|}(G(k,n);\mathbb Z)$.  The classes $\{\sigma_{\lambda}\}$ give an additive integral basis of  the cohomology ring of the Grassmannian.

Let $0< k_1 < \cdots < k_t < n$ be a sequence of positive integers. Let  $F(k_1, \dots, k_t; n)$ denote the partial flag variety which parameterizes partial flags $$F_{k_1} \subset F_{k_2} \subset \cdots \subset F_{k_t} \subset V_n,$$ where $\dim(F_{k_i}) = k_i$ for $1 \leq i \leq t$. When we would like to emphasize the ambient vector space, we will denote this partial flag variety by $F(k_1, \dots, k_t; V)$.

\subsection{Realizable Schubert classes}

Every Schubert class is represented by the corresponding Schubert variety, which is irreducible. The multiples of Schubert classes that are realizable over $\ZZ$ have been classified by Hong \cite{hong:rigidity} and Coskun and Robles \cite{coskunrobles}. Since we will frequently refer to these results in the proofs of our theorems, we state them here for the reader's convenience.

Express a partition $\lambda$ by grouping the equal parts together $(\mu_1^{i_1}, \dots, \mu_t^{i_t})$ so that 
$$\mu_1 > \mu_2 > \cdots > \mu_t$$ and precisely $i_j$ of the parts are equal to $\mu_j$ for $1 \leq j \leq t$. Recall that a Schubert class $\sigma_{\lambda}$ is called {\em multi rigid} if the only representatives of $m \sigma_{\lambda}$ for $m>0$ are the unions of $m$ Schubert varieties. In particular, if a Schubert class $\sigma_{\lambda}$ is multi rigid, then $m \sigma_{\lambda}$ is realizable over $\ZZ$ if and only if $m=1$.

\begin{theorem}\cite{hong:rigidity}
A Schubert class $\sigma_{\lambda}$ in $G(k,n)$ is multi rigid if and only if the following three conditions hold:
\begin{enumerate}
    \item $i_j \geq 2$ for $1 < j < t$ and
    \item $\mu_{j-1} \geq \mu_j + 2$ for $1 < j \leq t$ and 
    \item $i_1 \geq 2$ if $\mu_1 \not= n-k$ and $i_t \geq 2$ if $\mu_t \not= 0$.
\end{enumerate}
 \end{theorem}

\begin{theorem}\cite[Theorem 1.1]{coskunrobles}
    Let $\sigma_{\lambda}$ be a Schubert class in $G(k,n)$. Every positive integral multiple of $\sigma_{\lambda}$ can be represented by an irreducible subvariety if and only if $\sigma_{\lambda}$ is not multi rigid. 
\end{theorem}

Taken together these two theorems classify when multiples of Schubert classes in $G(k,n)$ are realizable over $\ZZ$.

\section{Realizable Classes}
In this section, we show that there are relations between the realization problem in different Grassmannians.  Fix integers $r\ge 1$,  $k\ge 1$ and $n\ge k$.

\begin{definition}
 Let  $\Lambda(r,k,n)$ be the set of partitions $\lambda = (\lambda_1,\ldots,\lambda_{k'})$ of length $\sum_i \lambda_i = r$ with $k'\le k$ parts satisfying
$$n-k \ge \lambda_1\ge \cdots \ge \lambda_{k'}\ge 1.$$
\end{definition}

We remark that we are taking the convention that our partitions have at most $k$ parts and do not have zeros.  Hence, we have inclusions 
$$\Lambda(r,k,n)\subset \Lambda(r,k+1,n+1).$$
Given any set $S$, let
$\Lambda(r,k,n;S)$ be the set of vectors  $(a_{\lambda})_{\lambda\in \Lambda(k,r)}$ of elements in $S$.

Observe that the cohomology classes of both codimension $r$ and dimension $r$ subvarieties of $G(k,n)$ are parameterized by $\Lambda(r,k,n; \mathbb Z_{\ge 0})$. Let $\sigma_\lambda$ denote the cohomology class of the Schubert cell corresponding to $\lambda$.  The class of a codimension $r$ subvariety is $$\sum_{\lambda \in \Lambda(r,k,n)} a_{\lambda} \sigma_{\lambda} \ \  \text{ with } a_{\lambda}\in \mathbb Z_{\ge 0}.$$
Let $\lambda^c$ denote the partition complementary to $\lambda$ in a $k(n-k)$ box; i.e., set $\lambda_i^c = n-k - \lambda_{k+1-i}$ for $1 \leq i \leq k$. We make the convention that $ \sigma^{\lambda}:=\sigma_{\lambda^c}$. Then we can write the class of a dimension $r$ subvariety of $G(k,n)$
as $$\sum_{\lambda \in \Lambda(r,k,n)} a_{\lambda^c} \sigma_{\lambda^c} = \sum_{\lambda \in \Lambda(r,k,n)} a^{\lambda} \sigma^{\lambda}.$$

\begin{definition} For $k,r\ge 1$ and $n\ge k$ set
\begin{align*}
I_{\mathbb Z}^r(k,n) &:= \left\{ \begin{array}{ll} (a_{\lambda})\in \Lambda(r,k,n;\mathbb Z_{\ge 0}) : & \text{the cohomology class}\\&  \sum a_{\lambda} \sigma_{\lambda} \in H^{2r}(G(k,n))\\ &\text{is  represented by an irreducible variety}\\& \text{of codimension $r$}\end{array}\right\}\\
I^{\mathbb Z}_r(k,n) &:= \left\{ \begin{array}{ll} (a^{\lambda})\in \Lambda(r,k,n;\mathbb Z_{\ge 0}) : & \text{the cohomology class } \\&\sum a^{\lambda} \sigma^{\lambda} \in H^{2(k(n-k)-r)}(G(k,n))\\ &\text{is  represented by an irreducible variety}\\ & \text{of dimension $r$}\end{array}\right\}\\
I_{\mathbb Q}^r(k,n) &:= \left\{ \begin{array}{ll} (a_{\lambda})\in \Lambda(r,k,n;\mathbb Q_{\ge 0}) : &\text{a positive  rational  multiple of the cohomology class } \\& \sum a_{\lambda} \sigma_{\lambda}\in H^{2r}(G(k,n)) \\ &\text{is represented by an irreducible variety}\\&\text{of codimension } r\end{array}\right\}\\
I^{\mathbb Q}_r(k,n) &:= \left\{ \begin{array}{ll} (a^{\lambda})\in \Lambda(r,k,n;\mathbb Q_{\ge 0}) : &\text{a positive  rational multiple of the cohomology class } \\& \sum a^{\lambda} \sigma^{\lambda}\in H^{2(k(n-k)-r)}(G(k,n)) \\ &\text{is represented by an irreducible variety}\\&\text{of dimension } r\end{array}\right\}
\end{align*}
We write $I^r(k,n)$ or $I_r(k,n)$ when we do not want to specify $\ZZ$ or $\QQ$.
\end{definition}

If  $(a_{\lambda})\in I_{\mathbb Z}^r(k,n)$, we say that \emph{$\sum_{\lambda \in \Lambda(r, k,n)} a_{\lambda} \sigma_{\lambda}$ is realizable over $\mathbb Z$}. When $(a_{\lambda})\in I_{\mathbb Q}^r(k,n)$, we say \emph{$\sum_{\lambda \in \Lambda(r, k,n)} a_{\lambda} \sigma_{\lambda}$ is realizable over $\mathbb Q$} (this terminology comes from \cite{hhmww}).

\subsection{Duality} The Grassmannian $G(k, V_n)$ is  isomorphic to the Grassmannian $G(n-k, V_n^*)$. By choosing an isomorphism between $V_n$ and $V_n^*$, there is a non-canonical  isomorphism $\phi$  between $G(k,n)$ and $G(n-k, n)$ that maps Schubert varieties  to Schubert varieties.  More precisely, given a partition $\lambda$, let $\lambda^T$ be its transpose. Then a Schubert variety with class $\sigma_{\lambda}$ maps to a Schubert variety with class $\sigma_{\lambda^T}$.
Hence, if $Y \subset G(k,n)$ has class $\sum_{\lambda\in \Lambda(k,n)} a_{\lambda} \sigma_\lambda$, then $\phi(Y)$ has class $\sum_{\lambda\in \Lambda(k,n)} a_{\lambda} \sigma_{\lambda^T}$. 

\begin{definition}
Let $a =(a_{\lambda}) \in \Lambda(r,k,n; S)$.  We define $a^T \in \Lambda(r,n-k,n; S)$ by the rule
$$ (a^T)_{\mu} : = a_{\mu^T} \text{ for } \mu\in \Lambda(r,n-k,n).$$ 
\end{definition}

We deduce the following proposition. 

\begin{proposition}\label{prop-duality}
If $a\in I^r(k,n)$ $($resp., $I_r(k,n)$$)$, then $a^T\in I^r(n-k,n)$ $($resp., $I_r(n-k,n)$$)$. 
\end{proposition}

\subsection{Inclusion of vector spaces}  The purpose of this subsection is to obtain inclusions between $I^r(k,n)$  (resp., $I_r(k,n)$) for different values of $k$ and $n$.

\begin{proposition}\label{prop-bumpkn}
We have the following inclusions:
\begin{enumerate}
    \item $I_r(k,n) \subseteq I_r(k, n+1)$.
    \item $I_r(k,n) \subseteq I_r(k+1, n+1)$.
    \item $I^r(k,n) \subseteq I^r(k+1, n+1)$.
    \item $I^r(k,n) \subseteq I^r(k, n+1)$.
\end{enumerate}
\end{proposition}

\begin{proof}
(1) We prove that $I_r(k,n) \subseteq I_r(k, n+1)$. Let $Y \subset G(k, V_n)$ be an irreducible subvariety of codimension $u$. Say $$[Y] = \sum_{\lambda\in \Lambda(u,k,n)}  a_{\lambda} \sigma_{\lambda}$$ in $H^{2u}(G(k,n))$.

For $s>0$, let $V_n \subset V_{n+s}$ be a linear inclusion. This inclusion induces an inclusion $G(k,V_n) \subset G(k, V_{n+s})$ and allows us to view $Y$ as a subvariety in $G(k, V_{n+s})$. Given a partition $\lambda$, let $\lambda + s^j$ be the partition obtained by adding $s$ to the first $j$ parts of  $\lambda$, in other words $$(\lambda + s^j)_i = \lambda_i + s \  \mbox{for} \  1 \leq i \leq j \quad  \mbox{and} \quad (\lambda + s^j)_i= \lambda_i \ \mbox{for} \ i>j.$$ Viewed in $G(k, V_{n+s})$, the variety $Y$ has class  $\sum a_{\lambda} \sigma_{\lambda+ s^k}$ in $H^{2u + 2ks} (G(k, n+s))$. 

In particular, if $Y \subset G(k,V_n)$ has dimension $r$ and class $\sum_{\lambda \in \Lambda(r, k,n)} a^{\lambda} \sigma^{\lambda}$, then $Y \subset G(k,V_{n+s})$ also has class $\sum_{\lambda \in \Lambda(r, k,n)} a^{\lambda} \sigma^{\lambda}$. We conclude that 
$I_r(k,n) \subseteq I_r(k, n+s)$ for any integer $s \geq 0$. Taking $s=1$, this proves (1).

 (2) We prove that  $I_r(k,n) \subseteq I_r(k+1, n+1)$. Let $t>0$ be an integer. Let $V_n \subset V_{n+t}$ be a linear inclusion and fix a subspace $V' \subset V_{n+t}$ of dimension $t$ such that $V' \cap V_n = 0$. Taking the spans $\overline{WV'}$ for $W \in G(k, V_n)$ induces an embedding of $G(k,n)$ in $G(k+t, n+t)$ and allows us to view a subvariety $Y$ of $G(k,n)$ as a subvariety of $G(k+t, n+t)$. Given a partition $\lambda \in \Lambda(u, k,n)$, let $(n-k)^j, \lambda$ denote the partition by appending $j$ parts equal to $n-k$ to $\lambda$, in other words, $$((n-k)^j, \lambda)_i = n-k \  \mbox{if} \ i \leq j \quad \mbox{and} \quad  ((n-k)^j, \lambda)_i = \lambda_{i-j} \ \mbox{if} \ i> j.$$ 
 So if $Y\subset G(k,n)$ is of dimension $u$ and has class $\sum a_{\lambda} \sigma_{\lambda}$ then viewed as a subvariety of $G(k+t, n+t)$ it has class  $\sum a_{\lambda} \sigma_{(n-k)^t, \lambda}$ in $H^{2u + 2t(n-k)} (G(k+t, n+t))$. 

In particular, if $Y \subset G(k,V_n)$ has dimension $r$ and class $\sum_{\lambda \in \Lambda(r, k,n)} a^{\lambda} \sigma^{\lambda}$, then $Y \subset G(k+t,V_{n+t})$ also has class $\sum_{\lambda \in \Lambda(r, k,n)} a^{\lambda} \sigma^{\lambda}$. We conclude that 
$I_r(k,n) \subseteq I_r(k+t, n+t)$ for any integer $t \geq 0$. This proves (2).

(3) We prove that $I^r(k,n) \subseteq I^r(k+1, n+1)$. Let $s, t >0$ be integers such that $t \leq n+s -k$.  For  an irreducible variety $Y\subset G(k,V_n)\subset G(k,V_{n+s})$  define the following incidence correspondence 
$$U:= \{ (W_1, W_2) \ | \ W_1 \in Y, \ W_2 \in G(k+t, V_{n+s}), W_1 \subset W_2\} \subset F(k, k+t; V_{n+s}).$$ The first projection $\pi_1: U \to G(k,V_{n+s})$ has image equal to $Y$. The fiber over a point $W_1 \in Y$ consists of all $(k+t)$-dimensional linear spaces containing $W_1$, which is the Grassmannian $G(t, V_{n+s}/W_1)$.  Since the incidence correspondence $U$ is a Grassmannian bundle over an irreducible variety, we conclude that $U$ is irreducible of dimension $\dim(Y) + t(n+s-k-t)$. 

Consider the  second projection $\pi_2: U \to G(k+t, V_{n+s})$. Let $W_2$ be a vector space in the image of $\pi_2$. The fiber of $\pi_2$ is the space of $W_1 \in Y$ that are contained in $W_2$. Now assume that $t \leq s$. Then the general $(k+t)$-plane  in $V_{n+s}$ containing a fixed $k$-plane in $V_n$ intersects $V_n$ in precisely that $k$-plane. To see this, complete a basis of $V_n$ to a basis of $V_{n+s}$ by adding vectors $e_{n+1}, \dots, e_{n+s}$. The span of $W_1$ with $e_{n+1}, \dots, e_{n+t}$ is such a $(k+t)$-plane. By semicontinuity, it follows that this property holds for the general $(k+t)$-plane. We conclude that $\pi_2$ is birational onto its image.  Hence $\pi_2(U)$ is irreducible of dimension $\dim(Y) + t(n+s-k-t).$

Note that $$k(n-k) - \mbox{codim}(Y) + t(n+s-k-t) = (k+t)(n+s-k-t) - \mbox{codim}(\pi_2(U)) - k(s-t).$$ In particular, if $s=t$, then the codimension of $Y$ in $G(k, V_n)$ is the same as the codimension of $\pi_2(U)$ in $G(k+t, V_{n+s})$.

We now compute the class of $\pi_2(U)$ in terms of the class of $Y$.
Let $[Y] = \sum a_{\lambda} \sigma_{\lambda}$ in $G(k, V_n)$. We claim that $[\pi_2(U)]= \sum a_{\lambda} \sigma_{\lambda + (s-t)^k}$ in $G(k+t, V_{n+s})$. In particular, when $t=s$, then the class of $Y$ and the class of $\pi_2(U)$ are denoted by the same partition in their respective Grassmannians.

To compute $[\pi_2(U)]$, we intersect with general Schubert varieties of complementary dimension. Suppose we intersect with a Schubert variety $\Sigma_{\mu}(G_{\bullet})$ with class $\sigma_{\mu}$ in $G(k+t, V_{n+s})$. First, in order to get a nonzero intersection, we  must have $\mu_1 = \cdots = \mu_t = n+s-k-t$. Otherwise, $\sum_{j>t} \mu_j > \dim(Y)$. Let $b_i = n+s-k-t+i - \mu_i$. Intersect the flag elements $G_{b_{t+1}} \subset \cdots \subset G_{b_{k+t}}$ with the vector space $V_n$. We get a partial flag in $V_n$ of dimensions $b_i - s$. Any $k$-plane parameterized by $Y$ must intersect $G_{b_i} \cap V_n$ in dimension at least $i-t$ for $t+1 \leq i \leq k+t$. This is a Schubert condition giving the class $\sigma_{\mu_{t+1}, \dots, \mu_{t+k}}$. Since $\sum_{j>t} \mu_j > \dim(Y)$, the corresponding Schubert variety does not intersect $Y$ and the intersection is $0$ as claimed. On the other hand, if $\mu_1 = \cdots = \mu_t = n+s-k-t$, the same argument shows that $\sum_{j>r} \mu_j = \dim(Y)$ and the number of intersection points is $a_{\lambda}$, where $\lambda$ is the dual of $\mu$. This proves the claim.

In particular, by taking $t=s=1$, we get that $I^r(k,n) \subseteq I^r(k+1, n+1)$. This proves (3).

 (4) We prove that $I^r(k,n) \subseteq I^r(k, n+1)$.  By Proposition \ref{prop-duality}, if $a \in I^r(k,n)$, then $a^T \in I^r(n-k, n)$. By the inclusion (3), $a^T \in I^r(n-k+1, n+1)$. By Proposition \ref{prop-duality}, $a\in I^r(k,n+1)$, as desired. 
\end{proof}

\subsection{The reverse inclusions}\label{sec:reverseinclusions} In this subsection, we prove the reverse inclusions, subject to certain restrictions on $k$ and $n$.

\begin{proposition}\label{prop-dumpdown} We have the following inclusions:
\begin{enumerate}
    \item If $n-k \geq r+1$, then $I_{\ZZ}^r(k, n+1) \subseteq I_{\ZZ}^r(k,n)$.
    \item If $n-k \geq r$, then $I_{\QQ}^r(k,n+1) \subseteq I_{\QQ}^r(k,n)$.
    \item If $k \geq r+2$, then $I_{\ZZ}^r(k, n+1) \subseteq I_{\ZZ}^r(k-1, n)$. 
    \item If $k \geq r+1$, then $I_{\QQ}^r(k, n+1) \subseteq I_{\QQ}^r(k-1, n)$.
    \item The only cohomology classes of codimension $r$ that are realizable over $\ZZ$ in $G(r+1, 2r+2)$ but are not realizable over $\ZZ$ in $G(r,2r)$ are  $$\{m \sigma_{r}, m \sigma_{1^r} \ | \  m >1\}.$$ 
\end{enumerate}
\end{proposition}
\begin{proof}
Parts (3) and (4) follow from parts (1) and (2), respectively, by duality. Let $Y$ be an irreducible subvariety of $G(k, n+1)$ with class $\sum a_{\lambda} \sigma_{\lambda}$. If $n + 1 - k \geq r+2$, then $\sigma_{2^k} \cdot [Y] \not=0$. So \cite[Theorem 8.1]{Debarre} implies the intersection of $Y$ with a general Schubert variety $\Sigma_{1^k}$ is irreducible. The Schubert variety $\Sigma_{1^k}$ is a Grassmannian $G(k,n)$ and the intersection $Y \cap \Sigma_{1^k}$ has class $\sum a_{\lambda} \sigma_{\lambda}$ in $G(k,n)$. This proves (1).  

If $n-k=r$ and $\sigma_{2^k} \cdot [Y]=0$, then $[Y]= a \sigma_{n-k}$. Since $\sigma_{n-k}$ is represented by a Schubert variety, we still have $I_{\QQ}^r (k, n+1) \subset I_{\QQ}^r(k, n)$. This proves (2). Observe that in this case, $a \sigma_{n-k}$ is represented by an irreducible subvariety in $G(k, n+1)$ \cite[Theorem 1.1]{coskunrobles}, but $a \sigma_{n-k}$ is not represented by an irreducible subvariety in $G(k,n)$ if $a > 1$ \cite{hong:rigidity}. Combining with duality, we conclude that the only cohomology classes that are realizable over $\ZZ$ in $G(r+1, 2r+2)$ but are not realizable over $\ZZ$ in $G(r,2r)$ are  $$\{m \sigma_{r}, m \sigma_{1^r} \ | \  m >1\}.$$ This concludes the proof of the proposition.
\end{proof}

\begin{corollary}\label{cor-equalitiesIJ} We have the following equalities.
\begin{enumerate}
\item If $n \geq k+r$, then $I_{\QQ}^r(k, n) = I_{\QQ}^r(k, k+r)$.
\item If $k \geq r$ and $n-k \geq r$, then $I_{\QQ}^r(k, n) = I_{\QQ}^r(r, 2r)$.
\item If $k \geq r+1$ and $n-k \geq r+1$, then $I_{\ZZ}^r(k, n) = I_{\ZZ}^r(r+1, 2r+2)$.

\end{enumerate}
\end{corollary}

\begin{proof}
By Proposition \ref{prop-bumpkn} (4), $I^r(k,n) \subseteq I^r(k, n+1)$. By Proposition \ref{prop-dumpdown} (2), if $n \geq k+r$, then $I_{\QQ}^r(k,n+1) \subseteq I_{\QQ}^r(k, n)$. This proves (1).

Let $e >0$ be an integer. By Proposition \ref{prop-bumpkn} (3), $I^r(r+e,2(r+e)) \subseteq I^r(r+e+j, 2r+2e+j)$ for every $j \geq 0$. By Proposition \ref{prop-bumpkn} (4), $I^r(r+e+j, 2r+2e+j) \subseteq I^r(r+e+j, n)$ for every $n \geq 2r+2e+j$. Setting $e=0$, we conclude that $$I^r(r, 2r) \subseteq I^r(k, n)$$ for every $k \geq r$ and $n \geq k+r$. Setting $e=1$, we deduce 
$$I^r(r+1, 2r+2) \subseteq I^r(k, n)$$ for every $k \geq r+1$ and $n \geq k+r+1$.

Conversely, suppose that $k=r+j$ and $n \geq 2r+j$ for some $j \geq 0$. Then Proposition \ref{prop-dumpdown} (4)  implies that $I^r_{\QQ}(k,n) \subseteq I_{\QQ}^r(r, n-j)$. Repeated applications of Proposition \ref{prop-dumpdown} (2) imply that $I_{\QQ}^r(r, n-j) \subseteq I_{\QQ}^r(r, 2r)$. We conclude that if $k \geq r$ and $n-k \geq r$, then $I_{\QQ}^r(k,n) = I_{\QQ}^r(r,2r)$. This proves (2).

Similarly, suppose that $k=r+1+j$ and $n \geq 2r+2+j$ for some $j \geq 0$. Then Proposition \ref{prop-dumpdown} (3) implies that $I^r_{\ZZ}(k,n) \subseteq I_{\ZZ}^r(r+1, n-j)$. Repeated applications of Proposition \ref{prop-dumpdown} (1) imply that $I_{\ZZ}^r(r+1, n-j) \subseteq I_{\ZZ}^r(r+1, 2r+2)$. We conclude that if $k \geq r+1$ and $n-k \geq r+1$, then $I_{\ZZ}^r(k,n) = I_{\ZZ}^r(r+1,2r+2)$. This proves (3).
\end{proof}

Analogous statements hold for dimension $r$ cohomology classes. 

\begin{proposition}\label{prop-dimdown} We have the following inclusions:
\begin{enumerate}
    \item If $n \ge r+k+1$, then $I_r(k, n+1) \subseteq I_r(k,n)$.
    \item If $n\ge r+k+1$ and $k\ge r+1$, then $I_r(k+1,n+1) \subseteq I_r(k,n)$.
    \item If $n \geq r+k$, then $I_r^{\QQ}(k, n+1) \subseteq I_r^{\QQ}(k,n)$.
    \item If $n \geq r+k$ and $k\geq r$, then $I_r^{\QQ}(k+1,n+1) \subseteq I_r^{\QQ}(k,n)$.
     \item The only cohomology classes of dimension $r$ that are realizable over $\ZZ$ in $G(r+1, 2r+2)$ but are not realizable over $\ZZ$ in $G(r,2r)$ are  $$\{m \sigma^{r}, m \sigma^{1^r} \ | \  m >1\}.$$ 
\end{enumerate}
\end{proposition}

\begin{proof}
    Interpret $G(k,n)$ as $\GG(k-1, n-1)$. An $r$-dimensional subvariety $Y$ of $G(k,n+1)$ induces an $(r+k-1)$-dimensional subvariety $Z$ of $\PP^{n}$ ruled by linear spaces $\PP^{k-1}$. If $r+k \leq n$, then the general point $p\in \PP^{n}$ is not contained in $Z$. Consequently, we can project $Z$ from $p$ to $\PP^{n-1}$ to obtain a subvariety $\pi_p(Z) \subset \PP^{n-1}$ ruled by linear spaces $\PP^{k-1}$. The variety $\pi_p(Z)$ induces a subvariety $Y'$ of $G(k, n)$. We claim that if $n>r+k$ and $p$ is general, then $Y$ and $Y'$ are birational. Take a general $\Lambda\cong \PP^{k-1}$ parameterized by $Y$. The space of $k$-dimensional projective linear spaces containing $\Lambda$ is $\PP^{n-k}$. If $n-k > r$, not every such $k$ plane can contain a $\PP^{k-1}$ parameterized by $Y$ since otherwise the dimension of $Y$ would be bigger than $r$. Finally, we note that if $[Y] = \sum a^{\lambda} \sigma^{\lambda}$, then $[Y'] = \sum a^{\lambda} \sigma^{\lambda}$ in their respective Grassmannians. This proves (1).

    Assume $n>r+k$ and $k>r$. By Proposition \ref{prop-duality}, if  $a \in I_r(k+1, n+1)$, then $a^T \in I_r(n-k, n+1)$. We have that $n > n-k +r$ since $k>r$. Hence part (1) of the proposition implies that $I_r(n-k, n+1) \subseteq I_r(n-k,n)$. By Proposition \ref{prop-duality}, if $a^T \in I_r(n-k,n)$, then $a \in I_r(k,n)$. Hence, $I_r(k+1, n+1)\subseteq I_r(k, n)$, proving (2).

    If $n-k=r$, it may happen that every $\PP^k$ containing a $\Lambda \cong \PP^{k-1}$ parameterized by $Y$ contains other $\PP^{k-1}$s parameterized by $Y$. Since $\dim(Y)=r$, for general choices, there can only be finitely many such $\PP^{k-1}$s. Hence, the map from $Y$ to $Y'$ is generically finite, say of degree $d$. In this case, we conclude that $[Y'] = \frac{1}{d} \sum a^{\lambda} \sigma^{\lambda}$ and $I_r^{\QQ}(k, n+1) \subseteq I_r^{\QQ}(k,n)$. This proves (3).  By duality we deduce (4).

    Let $Y \subset G(r+1, 2r+2)$. Note that if the map from $Y$ to $Y'$ is not birational, then any two linear spaces parameterized by $Y$ must span a $\PP^{r+1}$, hence must intersect in a $\PP^{r-1}$. A positive dimensional family of $\PP^{r}$s any two of which intersect in a $\PP^{r-1}$ must either contain a fixed $\PP^{r-1}$ or lie in a fixed $\PP^{r+1}$. In the first case, $Y$ has class $m \sigma^r$. Projecting to $G(r+1, 2r+1)$ gives an $m$ to 1 map to $\sigma^r$. Dually, in the second case, $Y$ has class $m \sigma^{1^r}$. In this case, $Y$ is a subvariety of $G(r+1, 2r+1)$. By duality, the dual of $Y$ is a subvariety of $G(r, 2r+1)$ with class $m \sigma^r$ and the discussion reduces to the previous case. The classes $m \sigma^r$ and $m \sigma^{1^r}$ are represented by irreducible subvarieties of $G(r+1, 2r+2)$ by \cite[Theorem 1.1]{coskunrobles} when $m \geq 1$. However, they  are not represented by an irreducible subvariety of $G(r,2r)$ if $m > 1$ by \cite{hong:rigidity}. This concludes the proof of the proposition.
\end{proof}

\begin{corollary}\label{cor-equalitydim} We have the following equalities:
\begin{enumerate}
    \item If $n \geq k+r$, then $I_r^{\QQ} (k,n) = I_r^{\QQ}(k, k+r)$.
     \item  Let $k \geq r$ and $n \geq r+k$, then $I_r^{\QQ}(k,n) = I_r^{\QQ}(r, 2r)$.
     \item Let $k \geq r+1$ and $n \geq r+k+1$, then $I_r^{\ZZ}(k,n) = I_r^{\ZZ}(r+1, 2r+2)$.
\end{enumerate}
  
\end{corollary}

\begin{proof}
By Proposition \ref{prop-bumpkn} (1), we have $I_r(k,n) \subseteq I_r(k, n+1)$. By Proposition \ref{prop-dimdown} (3), if $n \geq k+r$, $I_r^{\QQ}(k, n+1) \subseteq I_r^{\QQ}(k,n)$. This proves (1). 

   Fix an integer $e \geq 0$. By Proposition \ref{prop-bumpkn} (2), $I_r(r+e, 2r+2e+1) \subseteq I_r(r+e+j, 2r+2e+j)$ for every $j \geq 0$. By Proposition \ref{prop-bumpkn} (1), $I_r(r+e+j, 2r+2e+j) \subseteq I_r(r+e+j,n)$ for every $n \geq 2r+2e+j$. Taking $e=0$, we conclude that if $k \geq r$ and $n \geq r+k$, then $I_r(r, 2r) \subseteq I_r(k,n)$. Taking $e=1$, we conclude that if $k \geq r+1$ and $n \geq 2r+2$, then $I_r(r+1, 2r+2) \subseteq I_r(k,n)$.

    Conversely, suppose $k = r+1+j$ and $n \geq 2r+2+j$ for some $j \geq 0$. Then by Proposition \ref{prop-dimdown} (2), $I_r(k,n) \subseteq I_r(r+1, n-j)$. By repeated applications of Proposition \ref{prop-dimdown} (1), $I_r(r+1, n-j) \subseteq I_r(r+1, 2r+2)$. We conclude that if $k \geq r+1$ and $n\geq r+k+1$, then $I_r(k,n) = I_r(r+1, 2r+2)$. This proves (3)

    Moreover, if $k=r+j$ and $n \geq 2r+j$ for some $j \geq 0$, then $I_r^{\QQ}(k,n) \subseteq I_r^{\QQ}(r, n-j)$ by Proposition \ref{prop-dimdown} (4). By repeated applications of Proposition \ref{prop-dimdown} (3), $I_r(r, n-j) \subseteq I_r(r, 2r)$. We conclude that if $k \geq r$ and $n\geq r+k$, then $I_r^{\QQ}(k,n) = I_r^{\QQ}(r, 2r)$. This proves (2).
\end{proof}

\section{Construction of irreducible subvarieties of $G(k,n)$}

In this section, we construct irreducible subvarieties of $G(k,n)$ from irreducible subvarieties of products of smaller dimensional Grassmannians. This will be our main method of constructing irreducible subvarieties in $G(k,n)$.

The following proposition is well-known (see \cite{KleimanLandolfi} or \cite[Proposition 4.3]{coskun:rigid}).  We will sketch a proof for the convenience of the reader since it is central to our discussion and demonstrates the ideas of the general construction in a simpler setting. 

\begin{proposition}\label{prop-coneoversegre}
   The Pl\"{u}cker embedding of a Schubert variety $$\Sigma_{n^n, 0} \subset G(n+1, 2n+2)$$ is the cone (with vertex a point) over the Segre embedding of $\PP^n \times \PP^n$.  
\end{proposition}

\begin{proof}
    A Schubert variety $\Sigma_{n^n, 0}$ in $G(n+1, 2n+2)$ parameterizes $(n+1)$-dimensional subspaces of a $(2n+2)$-dimensional vector space $V$ that intersect a fixed $(n+1)$-dimensional space $\Lambda$ in a subspace of dimension at least $n$.  Fix an $(n+1)$-dimensional vector space $\Psi$ complementary to $\Lambda$. Given a vector $v \in \Psi$ and an $n$-dimensional vector space $U \subset \Lambda$, the span $W:=\overline{vU}$ is a point of $\Sigma_{n^n, 0}$. We thus get an injection  $f: \PP(\Psi) \times \PP(\Lambda)^* \to \Sigma_{n^n, 0}$. 

We claim that $$f^*(\OO_{G(n+1, 2n+2)}(1)) = \OO_{\PP(\Psi) \times \PP(\Lambda)^*}(1,1).$$ Fix an $(n+1)$-dimensional subspace $\Phi$ of $V$. The Schubert variety $\Sigma_{1, 0^n}(\Phi)$ parameterizing $(n+1)$-dimensional subspaces of $V$ that intersect $\Phi$ is the zero-locus of a  section of $\OO_{G(n+1, 2n+2)}(1)$. Choose $\Phi$ to intersect $\Lambda$ in a one-dimensional subspace spanned by $u$. By \cite[Theorem 3.39]{coskun:LR}, $f(\PP(\Psi) \times \PP(\Lambda)^*) \cap \Sigma_{1, 0^n}(\Phi)$ has two components, the first consisting of those $W$ that contain $u$ and the second consisting of those $W$ that lie in the $(2n+1)$-dimensional linear space $H$ spanned by $\Phi$ and $\Lambda$. In the first case, $U \subset \Lambda$ must contain the fixed vector $u$. Hence, the first component is  in the linear series of  $\pi_2^* (\OO_{\PP(\Lambda)^*}(1))$. In the second case, $v \in \Psi$ must lie in the hyperplane $H \cap \Psi$. Hence, the second component is in the linear series of $\pi_1^* (\OO_{\PP(\Psi)^*}(1))$. We conclude  that $f^*(\OO_{G(n+1, 2n+2)}(1)) = \OO_{\PP(\Psi) \times \PP(\Lambda)^*}(1,1)$. By choosing a basis $e_i$ with $1 \leq i \leq n+1$ for $\Lambda$ and $e_i$ with $n+2 \leq i \leq 2n+2$ for $\Psi$, it is easy to see that the span of  $f(\PP(\Psi) \times \PP(\Lambda)^*)$ has dimension $(n+1)^2 -1$ in $\PP(\bigwedge^n V)$. Hence, $f$ is the Segre embedding of $\PP(\Psi) \times \PP(\Lambda)^*$ under the complete linear system $|\OO_{\PP(\Psi) \times \PP(\Lambda)^*}(1,1)|$.

We now claim that $\Sigma_{n^n, 0}$ is the cone $\mathcal{C}$ over $f(\PP(\Psi) \times \PP(\Lambda)^*)$ with vertex $\Lambda$. The line joining $W$ and $\Lambda$ parameterizes vector spaces $W'$ that contain $U$ and are contained in the span $\overline{vU}$. Hence, the cone $\mathcal{C}$ is contained in $\Sigma_{n^n, 0}$.

Conversely, given any $(n+1)$-dimensional subspace $W' \not= \Lambda$ such that $[W']\in \Sigma_{n^n, 0}$, the projection of $W'$ from $\Lambda$ onto $\Psi$ determines a unique vector $v \in \Psi$ and $W' \cap \Lambda$ determines a unique $n$-dimensional subspace $U \subset \Lambda$. The vector space $W'$ contains $U$ and is contained in the vector space $\overline{v\Lambda}$, hence it lies on the line joining $W$ to $\Lambda$. We conclude that $\Sigma_{n^n, 0}$ is the cone with vertex $\Lambda$ over the image of $f$.
\end{proof}

For $1 \leq i \leq 2$, let $\pi_i: \PP^n \times \PP^n$ denote the $i$th projection. Let $$x : = [\pi_1^* (\OO_{\PP^n}(1)] \quad \mbox{and}  \quad y:= [\pi_2^* (\OO_{\PP^n}(1)].$$ Then the cohomology class of a $(2n-m)$-dimensional  subvariety $Z$ of $\PP^n \times \PP^n$ can be expressed as 
$$[Z]= \sum_{i=\max(0,m-n)}^{\min(m,n)} a_i x^i y^{m-i}.$$

\begin{theorem}\label{thm-mainconstruction}
    Let $Z$ be a $(2n-m)$-dimensional  irreducible subvariety of $\PP^n \times \PP^n$ with class $\sum_{i=\max(0,m-n)}^{\min(m,n)} a_i x^i y^{m-i}$. Then there exists an irreducible $(2n-m+1)$-dimensional subvariety of $G(n+1, 2n+2)$ with class $\sum_{i=\max(0,m-n)}^{\min(m,n)} a_i \sigma_{(n+1)^{m-i}, n^{n-m+i}, i}$.
\end{theorem}

\begin{proof}
We preserve the notation from the proof of Proposition \ref{prop-coneoversegre}.
 Given an irreducible variety $Z \in \PP(\Psi) \times \PP(\Lambda)^*$ of dimension $2n-m$, let $T$ be the cone over  $f(Z)$ with vertex $\Lambda$. Then $T$ is an irreducible $(2n-m+1)$-dimensional subvariety of $G(n+1, 2n+2)$. To compute the class of $T$, we can intersect $T$ by general Schubert varieties of complementary dimension $n^2+m$. First, observe that the coefficient of any class with $\lambda_n < n$ must be 0 since $T \subset \Sigma_{n^n, 0}$ and the dual of such a class  is disjoint from  $\Sigma_{n^n, 0}$. Hence, the class of $T$ is a linear combination of Schubert cycles of the form $\sigma_{(n+1)^{m-i}, n^{n-m+i}, i}$.

 Consider the cone over the image of a cycle with class $x^i y ^{m-i}$ under $f$. A representative of such a class is given by linear spaces $W$ that intersect $\Lambda$ in a fixed subspace $A$ of dimension $m-i$ and whose projection from $\Lambda$ lies in a fixed subspace $B$ of $\Psi$ of dimension $n+1-i$. Hence, such $W$ contain the $(m-i)$-dimensional subspace $A$, intersect $\Lambda$ in an $n$-dimensional subspace and are contained in the $(2n+2-i)$-dimensional subspace $\overline{B \Lambda}$. This is a Schubert variety  $\Sigma_{(n+1)^{m-i}, n^{n-m+i}, i}$. Hence, if the class of $Z$ is 
$$\sum_{i=\max(0,m-n)}^{\min(m,n)} a_i x^i y^{m-i},$$ then the class of $T$ is 
$$\sum_{i=\max(0,m-n)}^{\min(m,n)} a_i \sigma_{(n+1)^{m-i}, n^{n-m+i}, i}.$$ This concludes the proof of the theorem.
 \end{proof}

We can generalize Proposition \ref{prop-coneoversegre} in the following way.

\begin{proposition}\label{prop-bundleoversegre}
    Fix integers $0<i\leq n-k$ and $0 < j < k$. The Schubert variety $\Sigma_{i^j,0^{k-j}} \subset G(k, n)$ is the birational image of a $G(k-j, n-i-j)$-bundle over $G(j, n-k+j-i) \times G(k-j, k+i-j)$.
\end{proposition}

\begin{proof}
    The Schubert variety $\Sigma_{i^j,0^{k-j}} \subset G(k,n)$ parameterizes $k$-dimensional linear spaces $W$ of an $n$-dimensional vector space $V$ that intersect a fixed $(n-k+j - i)$-dimensional linear space $\Lambda$ in a subspace of dimension at least $j$. Fix a $(k+i-j)$-dimensional linear space $\Psi$ transverse to $\Lambda$. Consider the incidence correspondence
    $$U := \{ (W', W'', W) \ | \ W' \subset \Lambda, W'' \subset \Psi, W' \subset W \subset \overline{W'' \Lambda} \}$$ $$\subset G(j, n-k+j-i) \times G(k-j, k+i-j) \times G(k,n)$$ parameterizing a triple $(W', W'', W)$ of a $j$-dimensional subspace $W'$ of $\Lambda$, a $(k-j)$-dimensional subspace $W''$ of $\Psi$, and a $k$-dimensional subspace $W$ containing $W'$ and contained in the span of $W''$ and $\Lambda$. The third projection $\pi_3$ to the $W$ factor gives a map to $G(k,n)$. Since each such $W$ intersects $\Lambda$ in at least a $j$-dimensional subspace, the image of $\pi_3$ is contained in $\Sigma_{i^j, 0^{k-j}}$. Under $\pi_3$ the image of a fiber over a point $(W',W'')$ is the set of $k$-dimensional subspaces that contain $W'$ and are contained in $\overline{W''\Lambda}$. Observe that this is the Grassmannian $G(k-j,\overline{W''\Lambda}/W') \cong G(k-j, n-i-j)$. The map $\pi_3$ contracts the locus $$\{ (W', W'', W) \in U \ | \ \dim(W \cap \Lambda) \geq j+1 \}$$ to the singular locus of $\Sigma_{i^j,0^{k-j}}$.

    Conversely, given $W \in \Sigma_{i^j,0^{k-j}}$ such that $\dim(W \cap \Lambda)=j$, we can set $W \cap \Lambda = W'$  and $\overline{W \Lambda} \cap \Psi = W''$. Observe that $W''$ is a $(k-j)$-dimensional subspace of $\Psi$ and  $W \subset \overline{ W'W''}$. We conclude that $\pi_3$ is birational and its image is  $\Sigma_{i^j,0^{k-j}}$. This concludes the proof of the proposition.
\end{proof}

We can also generalize Theorem \ref{thm-mainconstruction}. By the K\"unneth decomposition, a codimension $r$ cohomology class in $G(j, n-k+j-i) \times G(k-j, k+i-j)$ can be expressed as $\sum a_{\lambda, \mu} \sigma_{\lambda} \otimes \sigma_{\mu}$, where the sum is over partitions $\lambda$ for $G(j, n-k+j-i)$ and $\mu$ for $G(k-j, k+i-j)$ such that $|\lambda| + |\mu|=r$.

\begin{theorem}\label{thm-mainconstruction2}
 Let $$Z \subset G(j, n-k+j-i) \times G(k-j, k+i-j)$$ be an $m$-dimensional irreducible subvariety  with class   $\sum a_{\lambda, \mu} \sigma_{\lambda} \otimes \sigma_{\mu}$. Then there exists an irreducible subvariety of $G(k,n)$ of dimension $m + (k-j) (n-i-k)$ with class
 $\sum a_{\lambda, \mu} \sigma_{\lambda + i^j, \mu}$.
\end{theorem}

\begin{proof}
    We preserve the notation from the proof of Proposition \ref{prop-bundleoversegre}. The projection $\pi_{1,2}$ from the incidence correspondence $U$ to $G(j, n-k+j-i) \times G(k-j, k+i-j)$ is a $G(k-j, n-i-j)$-bundle. Let $T=\pi_{1,2}^{-1}(Z)$. Then $T$ is an irreducible subvariety of $U$ of dimension $m+(k-j) (n-i-k)$. Since $\pi_3$ only contracts the locus of $W$ such that $\dim(W \cap \Lambda) \geq j+1$, the map $\pi_3$ restricted to $T$ is birational onto its image. Hence, $\pi_3(T)$ is an irreducible subvariety of $G(k,n)$ of dimension $m + (k-j) (n-i-k)$. 

    It remains to compute the class of $\pi_3(T)$. Let $S$ be a Schubert variety $S$ with class $\sigma_{\lambda} \otimes \sigma_{\mu}$ in $G(j, n-k+j-i) \times G(k-j, k+i-j)$. 
Given a flag $$F_1 \subset \cdots \subset F_{n-k+j-i} = \Lambda$$ in $\Lambda$ and a flag 
$$G_1 \subset \cdots \subset G_{k-i+j}= \Psi$$ in the complementary subspace $\Psi$, we get an induced flag on the ambient vector space $V$ by setting $$H_{\ell} = F_{\ell} \ \mbox{for} \ \ell \leq n-k+j-i \quad \mbox{and} \quad H_{\ell} = \overline{\Lambda G_{\ell - n+k-j+i}} \ \mbox{for} \ n-k+j-i < \ell \leq n.$$
The variety $\pi_3(\pi_{1,2}^{-1}(S))$ parameterizes $k$-dimensional subspaces of an $n$-dimensional vector space $V$ that intersect $F_{n-k-i+ \ell - \lambda_{\ell}}$ in a subspace of dimension at least $\ell$ for $1 \leq \ell \leq j$ and intersect $H_{n-k+\ell - \mu_{\ell- n+k-j+i}}$ in a subspace of dimension  at least $\ell$ for $j< \ell \leq k$.
This is precisely the Schubert variety with class $\sigma_{\lambda+ i^j, \mu}$ in $G(k,n)$. Hence, the class of $\pi_3(T)$ is $\sum a_{\lambda, \mu} \sigma_{\lambda + i^j, \mu}$. This concludes the proof of the theorem.
\end{proof}

By induction, Proposition \ref{prop-bundleoversegre} and Theorem \ref{thm-mainconstruction2} generalize to arbitrary Schubert varieties and subvarieties in products of more than two Grassmannians.

 Let $k < n$ be positive integers, and let $$n-k\geq i_1 > i_2 > \cdots > i_s >0$$ be a sequence of positive integers. Let $j_1, \dots, j_s$ be positive integers such that $\sum_{\ell=1}^s j_{\ell} < k$ and set $j_{s+1} = k - \sum_{\ell=1}^s j_{\ell}$. For $1 \leq t \leq s$, set $$a_{\ell} = n-k+\sum_{\ell=1}^t j_{\ell} - i_t.$$ Finally, set $a_0=0$ and $a_{s+1} = n$.

\begin{proposition}\label{prop-inductivestatement}
    The Schubert variety $\Sigma_{i_1^{j_1}, i_2^{j_2}, \dots, i_s^{j_s}, 0^{j_{s+1}}}$ is the birational image of a sequence of iterated  $G(j_{\ell}, a_{\ell -1} + j_{\ell} - \sum_{t=1}^{\ell-1} j_t)$-bundles for $2 \leq \ell \leq s+1$ over $\prod_{i=1}^{s+1} G(j_{\ell}, a_{\ell} - a_{\ell-1})$.
\end{proposition}

\begin{proof}
When $s=1$, this proposition reduced to Proposition \ref{prop-bundleoversegre}. 
    Fix a partial flag $$F_0 = 0 \subset F_{a_1} \subset \cdots \subset F_{a_s} \subset F_{a_{s+1}} = V,$$ where $F_{a_{\ell}}$ is an $a_{\ell}$-dimensional subspace of $V$. The Schubert variety $\Sigma_{i_1^{j_1}, i_2^{j_2}, \dots, i_s^{j_s}, 0^{j_{s+1}}}$ parameterizes $k$-dimensional subspaces $W$ such that $$\dim(W \cap F_{a_{t}}) \geq \sum_{\ell=1}^t j_{\ell} \ \ \mbox{for} \ 1 \leq t \leq s.$$ For each $1 \leq t \leq s$, fix a linear space $G_t \subset F_{a_{t+1}}$ transverse to $F_{a_t}$. Consider the incidence correspondence
$$U := \{ (W',Z_1, \dots, Z_s, W_{j_1}, W_{j_1+j_2}, \dots, W_k) \ |\ \dim(W') = j_1, \dim(Z_t) = j_{t+1},$$ $$\dim(W_{j_1 + \cdots + j_u}) = j_1 + \cdots + j_u,  \  W_{j_1} = W' \subset F_{a_1},\  Z_t \subset G_t, \ W_{j_1+ \cdots +j_t} \subset W_k \cap \overline{F_{a_{t-1}}Z_{t-1}}  \} $$ $$\subset \prod_{\ell=1}^{s+1} G(j_{\ell}, a_{\ell} - a_{\ell-1}) \times F(j_1, j_1+j_2, \dots, k;  n)$$ parameterizing tuples $(W', Z_1, \dots, Z_s, W_{j_1}, \dots, W_k)$, where
\begin{enumerate}
    \item $W'= W_{j_1}$ is a $j_1$-dimensional subspace of $F_{a_1}$,
    \item $Z_{\ell}$ is a $j_{\ell}$-dimensional subspace of $G_{\ell}$ for $1 \leq \ell \leq s$, 
    \item $W_{j_1} \subset W_{j_1+j_2} \subset \cdots \subset W_k$ is a partial flag in $V$, where the $u$th element of the partial flag has dimension $\sum_{\ell=1}^u j_{\ell}$  and is contained in the span of $F_{a_{u-1}}$ and $Z_{u-1}$.
\end{enumerate}

The projection to the last factor $W_k$ gives a map from $U$ to $G(k,n)$. By construction, $$\dim(W \cap F_{a_t}) \geq \sum_{\ell=1}^t j_{\ell} \ \  \mbox{for} \ 1 \leq t \leq s.$$ Hence, the image is contained in the Schubert variety $\Sigma_{i_1^{j_1}, i_2^{j_2}, \dots, i_s^{j_s}, 0^{j_{s+1}}}$. Conversely, given a point $W \in \Sigma_{i_1^{j_1}, i_2^{j_2}, \dots, i_s^{j_s}, 0^{j_{s+1}}}$ such that $\dim (W \cap F_{a_t}) = \sum_{\ell=1}^t j_{\ell}$ for every $1 \leq t \leq s$, we can recover $W_{j_1 + \cdots + j_t}$ as $W \cap F_{a_t}$. Projecting $W \cap F_{a_t}$ from $F_{a_{t-1}}$ to $G_{t-1}$ uniquely determines the linear space $Z_t$. Hence, the last projection maps $U$ birationally onto $\Sigma_{i_1^{j_1}, i_2^{j_2}, \dots, i_s^{j_s}, 0^{j_{s+1}}}$. Successively forgetting the elements of the partial flag realizes $U$ as an iterated Grasmannian bundle over $\prod_{i=1}^{s+1} G(j_{\ell}, a_{\ell} - a_{\ell-1})$. If we have chosen $W_{j_1 + \cdots + j_{u-1}}$, the choice of $W_{j_1+ \cdots + j_u}$ is a choice of $(j_1 + \cdots + j_u)$-dimensional linear space containing $W_{j_1 + \cdots + j_{u-1}}$ and contained in $\overline{F_{a_{u-1}}Z_{u-1}}$. This choice is parameterized by the Grassmannian $G(j_u, a_{u-1} + j_u - \sum_{\ell=1}^{u-1} j_{\ell})$. This concludes the proof of the proposition.
\end{proof}

\begin{theorem}\label{thm-mainconstruction2ind}
    Let $Z$ be an $m$-dimensional irreducible subvariety of $\prod_{\ell=1}^{s+1}G(j_{\ell}, a_{\ell} - a_{\ell-1})$ with class $\sum a_{\lambda^1, \dots, \lambda^{s+1}}  \sigma_{\lambda^1}\otimes \cdots \otimes \sigma_{\lambda^{s+1}}$. Then there exists an $(m + \sum_{\ell=2}^{s+1} j_{\ell} (a_{\ell-1} - \sum_{t=1}^{\ell-1} j_t))$-dimensional irreducible subvariety of $G(k,n)$ with class
    $$\sum a_{\lambda^1, \dots, \lambda^{s+1}} \sigma_{\lambda^1 + i_1^{j_1}, \lambda^2 + i_2^{j_2}, \dots, \lambda^s + i_s^{j_s}, \lambda^{s+1}}.$$
\end{theorem}

\begin{proof}
 We preserve the notation from Proposition \ref{prop-inductivestatement}. Let $\pi$ denote the projection from $U$ to the first $s+1$ factors. Let $\phi$ denote the projection to $G(k,n)$. Given $Z \in \prod_{\ell=1}^{s+1}G(j_{\ell}, a_{\ell} - a_{\ell-1})$, consider $T=\phi(\pi^{-1}(Z))$. The map $\phi$ is birational onto its image when restricted to  $\pi^{-1}(Z)$. Hence, $T$ is an irreducible variety of the stated dimension. The class of $T$ can be computed as in the proof of Theorem \ref{thm-mainconstruction2}. 
\end{proof}

\section{Dimension 2 and Codimension 2 realizable classes}
In this section, we determine all the cohomology classes in $G(k,n)$ of dimension or codimension 2 that can be represented by irreducible subvarieties.

\begin{theorem}\label{thm-dimension2}
Let $2 \leq k \leq n-k$ and let $$\nu= a \sigma_{(n-k)^{k-1}, n-k-2} + b \sigma_{(n-k)^{k-2}, (n-k-1)^2}= a\sigma^2 + b \sigma^{1,1}$$ be a non-zero, integral cohomology class.
\begin{enumerate}
    \item When $k >2$, the class $\nu$ can be represented by an irreducible surface if and only if $a, b \geq 0$.
    \item When $k=2$ and $n>4$, the class $\nu$ can be represented by an irreducible surface if and only if $a > 0$ and $b \geq 0$ or $(a,b)=(0,1)$.
\item When $(k,n)=(2,4)$, the class $\nu$ can be represented by an irreducible surface if and only if $a, b > 0$ or $(a,b)=(1,0)$ or $(a,b)=(0,1)$.
    \end{enumerate}
\end{theorem}

\begin{proof}
The effective cone of cycles in a Grassmannian is spanned by classes of Schubert cycles \cite[Proposition 2.20]{coskun:IMPANGA}. Hence, the class of any effective cycle is a nonnegative integral combination of Schubert cycles.

First, consider the case of $G(2,4)$. By Proposition \ref{prop-coneoversegre}, a Schubert variety $\Sigma_{1,0}$ in $G(2,4)$ is the cone over $\PP^1 \times \PP^1$. By Bertini's Theorem, there exists irreducible curves of type $(a,b)$ on $\PP^1 \times \PP^1$, provided that either $a,b >0$; or $a=1$ and $b=0$; or $a=0$ and $b=1$. Hence, the claimed classes can be represented by irreducible surfaces by Theorem \ref{thm-mainconstruction}. Conversely, the Schubert classes $\sigma_{2,0}$ and $\sigma_{1,1}$ are multi rigid in $G(2,4)$ \cite{hong:rigidity}, hence $m \sigma_{2,0}$ or $m \sigma_{1,1}$ can be represented by an irreducible subvariety if and only if $m=1$. This concludes the discussion in the case of $G(2,4)$.

If $2 \leq k \leq n-k$ and $(k,n) \not= (2,4)$, by repeated applications of Proposition \ref{prop-bumpkn} (1) and (2), we have $I_2^{\ZZ}(2,4) \subseteq I_2^{\ZZ}(k,n)$. 
When $k=2$ and $n>4$, then the class $m \sigma_{(n-k)^{k-1}, n-k-2} = m \sigma^2$ is representable by an irreducible surface by \cite[Theorem 1.1]{coskunrobles}. On the other hand, the class $m \sigma_{(n-k)^{k-2}, (n-k-1)^2}= m \sigma^{1,1}$ is  representable by an irreducible surface if and only if $m=1$ by \cite{hong:rigidity}. This concludes the discussion when $k=2$.

When $k>2$, the classes $m \sigma_{(n-k)^{k-1}, n-k-2} = m \sigma^2$ and $m \sigma_{(n-k)^{k-2}, (n-k-1)^2}= m \sigma^{1,1}$ are representable by irreducible subvarieties by  \cite[Theorem 1.1]{coskunrobles}. Hence, every integral cohomology class $a \sigma^2 + b \sigma^{1,1}$ with $a, b \geq 0$ can be represented by an irreducible surface in $G(k,n)$. This concludes the proof. 
\end{proof}

A similar argument classifies the codimension 2 classes that can be represented by irreducible subvarieties. 

\begin{theorem}\label{thm-codimension2}
Let $2 \leq k \leq n-k$ and let $\nu= a \sigma_{2} + b \sigma_{1,1}$ be a non-zero, integral cohomology class.
\begin{enumerate}
    \item When $k >2$, the class $\nu$ can be represented by an irreducible subvariety if and only if $a, b \geq 0$.
    \item When $k=2$ and $n>4$, the class $\nu$ can be represented by an irreducible subvariety if and only if $a > 0$ and $b \geq 0$ or $(a,b)=(0,1)$.
\item When $(k,n)=(2,4)$, the class $\nu$ can be represented by an irreducible subvariety if and only if $a, b > 0$ or $(a,b)=(1,0)$ or $(a,b)=(0,1)$.
    \end{enumerate}
\end{theorem}

\begin{proof}
Since $G(2,4)$ has dimension 4, we have $I_2(2,4) = I^2(2,4)$. Hence, when $(k,n)=(2,4)$, Theorem \ref{thm-dimension2} and Theorem \ref{thm-codimension2} coincide. 

If $2 \leq k \leq n-k$ and $(k,n) \not= (2,4)$, by repeated applications of Proposition \ref{prop-bumpkn} (3) and (4), we have $I^2_{\ZZ}(2,4) \subseteq I^2_{\ZZ}(k,n)$. Moreover, if $k>2$, the classes $m \sigma_{2}$ and $m \sigma_{1,1}$ are representable by irreducible subvarieties by  \cite[Theorem 1.1]{coskunrobles} for $m \geq 1$. Hence, every effective, integral cohomology class $a \sigma^2 + b \sigma^{1,1}$ with $a, b \geq 0$ can be represented by an irreducible subvariety in $G(k,n)$. When $k=2$, the classes $m \sigma_{2}$ are representable by irreducible subvarieties by  \cite[Theorem 1.1]{coskunrobles} for $m \geq 1$. On the other hand, by \cite{hong:rigidity}, $m \sigma_{1,1}$ is representable by an irreducible subvariety if and only if $m=1$. This concludes the proof. 
\end{proof}

\section{Dimension 3 realizable classes}
In this section, we classify three-dimensional realizable classes in $G(k,n)$.

\begin{theorem}\label{thm-dimension3}
 Let $3 \leq k \leq n-k$ and let $$\nu=  a \sigma^3 +  b \sigma^{2,1} + c \sigma^{1,1,1}$$ be a nonzero integral cohomology class.
 \begin{enumerate}
     \item If $k >3$, then $\nu$ is realizable over $\ZZ$ if and only if $a,b,c \geq 0$ and $b^2 \geq ac$.
     \item If $k=3$ and $n>6$, then $\nu$ is realizable over $\ZZ$ if and only if $a,b,c \geq 0$, $b^2 \geq ac$ and one of the following three conditions holds: 
     \begin{enumerate}
         \item[(i)] $b>0$; or
         \item[(ii)] $b=c=0$; or 
         \item[(iii)] $a=b=0$ and $c=1$.
     \end{enumerate}
     \item If $k=3$ and $n=6$, then $\nu$ is realizable over $\ZZ$ if and only if $a,b,c \geq 0$, $b^2 \geq ac$ and one of the following three conditions holds: 
     \begin{enumerate}
         \item[(i)] $b>0$; or 
         \item[(ii)] $(a,b,c)=(1,0,0)$; or 
         \item[(iii)] $(a,b,c) = (0,0,1)$.
     \end{enumerate}
    \end{enumerate}
 \end{theorem}

\begin{proof}
An effective class is a nonnegative linear combinations of Schubert classes \cite[Proposition 2.20]{coskun:IMPANGA} , hence $a,b,c \geq 0$.

    First consider the case of $G(3,6)$. The Schubert variety $\Sigma_{2,2,0}$ is the cone over $\PP^2 \times \PP^2$ by Proposition \ref{prop-coneoversegre}. By \cite[Theorem 1]{huh2} a surface class $a x^2 + bxy + cy^2$ in $\PP^2 \times \PP^2$, where $x$ and $y$ are the pullbacks of the hyperplane classes from the two factors, is realizable over $\ZZ$ if and only if $a,b,c \geq 0$ and one of the following holds: (i) $b>0$ or (ii) $(a,b,c)=(1,0,0)$ or (iii) $(a,b,c)=(0,0,1)$. By Theorem \ref{thm-mainconstruction}, $a \sigma^3 + b \sigma^{2,1} + c \sigma^{1,1,1}$ is realizable over $\ZZ$ if $a,b,c \geq 0$, $b^2 \geq ac$ and one of the following three conditions holds: (i) $b>0$ or (ii) $a=1$, $b=c=0$ or (iii) $a=b=0$ and $c=1$. This proves that the classes listed in (3) are realizable over $\ZZ$.

    When $k \geq 3$ and $2k \leq n$, by repeated applications of Proposition \ref{prop-bumpkn} (1) and (2), $$I_3^{\ZZ} (3, 6) \subseteq I_3^{\ZZ}(k,n).$$ Assume that $k=3$. If $n=6$, then $m \sigma^3$ or $m \sigma^{1,1,1}$ is realizable over $\ZZ$ if and only if $m=1$ by \cite{hong:rigidity}. If $n >6$, then $m \sigma^3$ is realizable over $\ZZ$ for any $m \geq 1$ by \cite[Theorem 1.1]{coskunrobles}. On the other hand, $m \sigma^{1,1,1}$ is realizable over $\ZZ$ if and only if $m=1$ by \cite{hong:rigidity}. Finally, if $k>3$, then $m \sigma^3$ or $m \sigma^{1,1,1}$ is realizable over $\ZZ$ for any $m\geq 1$ by \cite[Theorem 1.1]{coskunrobles}. This proves that the classes listed in the theorem are realizable over $\ZZ$.

    There remains to show that no other three-dimensional classes can be represented by an irreducible subvariety of $G(k,n)$. We claim it suffices to show that if $\nu$ is realizable over $\ZZ$, then $b^2 \geq ac$. If $b>0$, we already proved that the classes satisfying $b^2 \geq ac$ are realizable over $\ZZ$. On the other hand, if $b=0$, then either $a=0$ or $c=0$ since $a, c\geq 0$. The previous paragraph classifies the classes of the form $m \sigma^3$ or $m \sigma^{1,1,1}$ that are realizable over $\ZZ$.

    By Corollary \ref{cor-equalitydim} (2), if $k\geq 3$ and $n \geq 2k$, then $I_3^{\QQ}(k,n) = I_3^{\QQ}(3,6)$. Hence, it suffices to show $b^2 \geq ac$ when $k=3$ and $n=6$.

    Let $V$ be a $6$-dimensional vector space. Fix a $5$-dimensional subspace $A$ and a $4$-dimensional subspace $B$ of $V$.  Consider the birational map 
    $$f: G(2,A) \times \PP(B) \dashrightarrow G(3,6)$$ obtained by taking a two-dimensional subspace $W' \subset A$ and a one-dimensional subspace $p \subset B$ and mapping it to their span $W=\overline{pW'}$. This is a birational map between $G(2,5) \times \PP^3$ and $G(3,6)$. Given an irreducible, three-dimensional subvariety $Z \subset G(3,6)$, we can pull back a general translate of $Z$ to $G(2,5) \times \PP^3$ to obtain an irreducible  threefold in $G(2,5) \times \PP^3$. 

Let $x$ denote the hyperplane class on $\mathbb P(B)$.  We claim that the classes on $G(3,6)$ pullback as follows:
\begin{align*}
& f^* \sigma^3 \ \ \ \ = \sigma^0 + \sigma^1 x + \sigma^2 x^2 + \sigma^3 x^3 \\
&f^* \sigma^{2,1} \ \ = \sigma^1 x + \sigma^2 x^2 + \sigma^{1,1} x^2 + \sigma^{2,1} x^3  \\
&f^* \sigma^{1,1,1} = \sigma^{1,1}x^2
    \end{align*}

Here and in subsequent proofs these pullbacks can be computed by elementary Schubert calculus. To compute the coefficient of $\sigma^{i,j} x^k$ in a class $f^* \nu$, we can intersect $f^* \nu$ with the dual Schubert class $\sigma^{3-j, 3-i} x^{3-k}$. This amounts to computing the intersection product $\sigma_{i,j} \cdot \sigma_{3-k} \cdot \nu$ in $H^*(G(3,6), \ZZ)$. These are readily computed using any Littlewood-Richardson rule (see, for example, \cite[Theorem 3.39]{coskun:LR}).   
    
Intersecting $Z$ with the nef class $x$, we obtain a class which is a limit of irreducible surfaces. This has class 
$$[Z]x= a(\sigma^0 x + \sigma^1 x^2 + \sigma^2 x^3 ) + b (\sigma^1 x^2 + \sigma^2 x^3 + \sigma^{1,1} x^3) + c \sigma^{1,1}x^3.$$
We can write the intersection matrix with respect to $x$ and $\sigma_1$, to obtain 
$$\left( \begin{array}{cc} 
a & a+b \\
a+b & a+ 2b+c
\end{array}\right). $$
By the Hodge index theorem, the determinant of this matrix has to be nonpositive. We conclude that 
$$a(a+2b+c) - (a+b)^2 = ac-b^2 \leq 0.$$ This proves that $b^2 \geq ac$, completing the proof of the theorem.
\end{proof}

For completeness, we discuss dimension 3 realizable classes in the Grassmannians $G(2,n)$. A three-dimensional class in $G(2,4)$ is a divisor with class $m \sigma_1 = m \sigma^{2,1}$. If $m \geq 1$, such a class is represented by an irreducible hypersurface. Hence, we can assume that $n \geq 5$.

\begin{theorem}\label{thm-dimension3k2}
    Let $\nu = a \sigma^3 + b \sigma^{2,1}$ be a nonzero, integral cohomology class. 
    \begin{enumerate}
        \item If $n > 5$, then $\nu$ is realizable over $\ZZ$ in $G(2, n)$ if and only if $a, b \geq 0$.
        \item The class $\nu$ is realizable over $\ZZ$ in $G(2, n)$ if and only if either $a \geq 0$ and $b >0$; or $(a,b) = (1,0)$.
    \end{enumerate}
\end{theorem}

\begin{proof}
If $\nu$ is effective, then $a, b \geq 0$.

   The Schubert variety $\Sigma_2 \in G(2,5)$ is a cone over the Segre embedding of $\PP^1 \times \PP^2$ by Proposition \ref{prop-coneoversegre}.  Let $x,y$ denote the hyperplane class of $\mathbb P^1$ and $\mathbb P^2$, respectively.  Given a surface $S\subset \mathbb P^1\times \mathbb P^2$ with class $a x + by$,  the cone $T$ over $S$ is a threefold with class $a \sigma^3 + b \sigma^{2,1}$ in $G(2,5)$ by Theorem \ref{thm-mainconstruction}. Since the class $ax + by$ is very ample on $\PP^1 \times \PP^2$ when $a, b >0$, we conclude that it can be represented by an irreducible surface. Hence, $a \sigma^3 + b \sigma^{2,1}$ is realizable over $\ZZ$ in $G(2,5)$ if $a,b >0$. The class $b \sigma^{2,1}$ is realizable over $\ZZ$ for $b >0$ by \cite[Theorem 1.1]{coskunrobles}. The class $a \sigma^3$ is realizable over $\ZZ$ if and only if $a=1$ by \cite{hong:rigidity}. This proves (2).

   By Proposition \ref{prop-bumpkn} (1), $I_3^{\ZZ}(2,5) \subseteq I_3^{\ZZ} (2,n)$ for $n > 5$. The classes $m \sigma^3$ and $m \sigma^{2,1}$ are realizable over $\ZZ$ for every $m > 0$ by \cite[Theorem 1.1]{coskunrobles}. Hence, every nonzero effective class $a \sigma^3 + b \sigma^{2,1}$ is realizable over $\ZZ$ in $G(2,n)$ when $n>5$. This proves (1).
\end{proof}

\section{Codimension 3 realizable classes}

In this section, we classify the codimension 3 classes in $G(k,n)$ that are realizable over $\ZZ$.

\begin{theorem}\label{thm-codimension3}
Let $3 \leq k \leq n-k$. Let $\nu = a \sigma_3 + b \sigma_{2,1} + c \sigma_{1,1,1}$ be a nonzero, integral cohomology class in $G(k,n)$.
\begin{enumerate}
    \item If $k > 3$, then $\nu$ is realizable over $\ZZ$ in $G(k,n)$ if and only if $a,b,c \geq 0$ and $b^2 \geq ac$.
    \item If $k=3$ and $n>6$, then $\nu$ is realizable over $\ZZ$ in $G(k,n)$ if and only if $a,b,c \geq 0$, $b^2 \geq ac$ and one of the following three conditions holds:
    \begin{enumerate}
        \item[(i)] $b>0$; or
        \item[(ii)] $b=c=0$ and $a>0$; or
        \item[(iii)] $(a,b,c)=(0,0,1)$.
    \end{enumerate}
    \item The class $\nu$ is realizable over $\ZZ$ in $G(3,6)$ if and only if $a,b,c \geq 0$, $b^2 \geq ac$ and one of the following three conditions holds:
    \begin{enumerate}
        \item[(i)] $b>0$; or
        \item[(ii)] $(a,b,c)=(1,0,0)$; or
        \item[(iii)] $(a,b,c)=(0,0,1)$.
    \end{enumerate}
\end{enumerate}
\end{theorem}

\begin{proof}
  If $\nu$ is an effective class, then $a,b,c \geq 0$ by \cite[Proposition 2.20]{coskun:IMPANGA}. Hence, below we will assume that $a,b,c$ are nonnegative.

  We first discuss the case of $G(3,6)$. Fix two transverse $3$-dimensional linear spaces $A$ and $B$. Consider the incidence correspondence
    $$ U := \{ (p, l, W) \ |\ p\in A, l \subset B, p \in W \subset \overline{lA} \} \subset \PP(A) \times \PP(B)^* \times G(3,6)$$
    parameterizing triples $(p, l, W)$, where $p$ is a one-dimensional subspace of $A$ contained in $W$ and  $l$ is a two-dimensional subspace of $B$ such that $W$ is contained in the span of $l$ and $A$. By Proposition \ref{prop-bundleoversegre}, $U$ is a $G(2,4)$-bundle over $\PP^2 \times (\PP^2)^*$. Let $x$ and $y$ denote the pullbacks of $\OO_{\PP^2}(1)$ and $\OO_{(\PP^2)^*}(1)$, respectively. By Theorem \ref{thm-mainconstruction2}, given a surface with class $a y^2 + b xy + cy^2$ in $\PP^2 \times (\PP^2)^*$, the inverse image under $\pi_{1,2}$ followed by the projection by $\pi_3$ is an irreducible sixfold in $G(3,6)$ with class
    $a \sigma_3 + b \sigma_{2,1} + c \sigma_{1,1,1}$.  Since the class $a y^2 + b xy + cy^2$ is realizable over $\ZZ$ in $\PP^2 \times (\PP^2)^*$ when $b>0$ and $b^2 \geq ac$ by \cite[Theorem 1]{huh2}, we conclude that the classes listed in (3) are realizable over $\ZZ$.

    By repeated applications of Proposition \ref{prop-bumpkn} (3) and (4), if $k \geq 3$ and $n \geq k$, we have $I_{\ZZ}^3(3,6) \subseteq I_{\ZZ}^3(k,n)$. Hence, the classes $a y^2 + b xy + cy^2$ with $a,c \geq 0$, $b>0$ and $b^2 \geq ac$ are realizable over $\ZZ$ in any $G(k,n)$ with $3 \leq k \leq n-k$. By \cite[Theorem 1.1]{coskunrobles}, $a \sigma_3$ is represented by an irreducible subvariety for $a >0$ in $G(k,n)$ if $k \geq 3$ and $n>6$ and $c \sigma_{1,1,1}$ is represented by an irreducible subvariety for $c >0$ in $G(k,n)$ if $k > 3$ and $n\geq 2k$. Hence, the classes listed in the theorem are realizable over $\ZZ$.

    To conclude the proof of the theorem, we need to show that the classes listed in the theorem are the only classes that are realizable over $\ZZ$. We claim that for this it suffices to show that if $\nu$ is realizable over $\ZZ$, then $b^2 \geq ac$. We have already shown that if $b>0$ and $b^2 \geq ac$, then $\nu$ is realizable over $\ZZ$. On the other hand, if $b=0$ and $b^2 \geq ac$, then either $a=0$ or $c=0$. If $k=3$, then $c \sigma_{1,1,1}$ is realizable over $\ZZ$ if and only if $c=1$ by \cite{hong:rigidity}. If $k=3$ and $n=6$, then $a \sigma_3$ is realizable over $\ZZ$ if and only if $a=1$ by \cite{hong:rigidity}. This completes the classification of the realizable classes when $b=0$ subject to proving the inequality $b^2 \geq ac$.

    By Corollary \ref{cor-equalitiesIJ} (2), $I_{\QQ}^3(k,n) = I_{\QQ}^3(3,6)$ when $3 \leq k \leq n-k$. Hence, it suffices to prove the inequality $b^2 \geq ac$ for $G(3,6)$. Let $Y$ be an irreducible sixfold in $G(3,6)$ with class $a \sigma_3 + b \sigma_{2,1} + c \sigma_{1,1,1}$. We can linearly embed $V_6 \subset V_8$. Fix a two-dimensional linear space $\Lambda \subset V_8$ transverse to $V_6$. By sending $W \in Y$ to $\overline{\Lambda W}$, we get an embedding of  $Y$ in $G(5,8)$ with class 
    $$[Y]= a\sigma_{3,3,3} + b \sigma_{3,3,2,1} + c \sigma_{3,3,1,1,1}.$$
    We now consider the birational map $$f : G(2,5) \times G(3,6) \dashrightarrow G(5,8)$$ obtained by fixing two transverse linear spaces $\Omega, \Gamma \subset V_8$ of dimensions 5 and 6, respectively, and sending $W \in G(2,\Omega)$ and $W' \in G(3, \Gamma)$ to $\overline{WW'}$. Pulling back $Y$ by $f$, we obtain a sixfold in $G(2,5) \times G(3,6)$. We now apply the Hodge index theorem on the class $\sigma_{2,2}\otimes 1 \cdot [f^{-1}(Y)]$. We have the following 
    \begin{align*}
        & \sigma_{2,2}\otimes 1 \cdot f^* \sigma_{3,3,3} \ \ \ \ = \sigma_{2,2} \otimes \sigma_{3,3,3} + \sigma_{3,2}\otimes \sigma_{3,3,2} + \sigma_{3,3} \otimes \sigma_{3,2,2} \\
    & \sigma_{2,2}\otimes 1 \cdot  f^* \sigma_{3,3,2,1} \ \ = \sigma_{3,2}\otimes \sigma_{3,3,2} + \sigma_{3,3} \otimes \sigma_{3,3,1} + \sigma_{3,3} \otimes \sigma_{3,3,2}\\
    & \sigma_{2,2}\otimes 1 \cdot f^* \sigma_{3,3,1,1,1} = \sigma_{3,3} \otimes \sigma_{3,3,1}
    \end{align*}
Hence, $$\sigma_{2,2}\otimes 1 \cdot [f^{-1}(Y)] = a \sigma_{2,2}\otimes \sigma_{3,3,3} + (a+b)\sigma_{3,2}\otimes \sigma_{3,3,2} + (a+b)\sigma_{3,3} \otimes \sigma_{3,2,2} + (b+c) \sigma_{3,3} \otimes \sigma_{3,3,1}.$$
Calculating the intersection matrix for $\sigma_1 \otimes 1$ and  $1 \otimes \sigma_1$, we obtain
$$\left(\begin{array}{
cc}
    a & a+b \\
   a+b  & a+2b+c
\end{array}\right).$$
By the Hodge index theorem, the determinant of this matrix must be non-positive:
$$a(a+2b+c)-(a+b)^2 = ac-b^2 \leq 0.$$ Hence, we conclude $b^2 \geq ac$ as desired. This completes the proof. 

\end{proof}

For completeness, we describe the realizable codimension $3$ classes in $G(2,n)$. In $G(2,4)$, a codimension 3 class is a curve,  and every cohomology class $m \sigma_{2,1}$ with $m > 0$ can be represented by an irreducible curve. In $G(2,5)$, codimension 3 and dimension 3 agree, so Theorem \ref{thm-dimension3k2} (2) classifies the classes that are realizable over $\ZZ$. We may, therefore, assume that $n \geq 6$.
 
\begin{theorem}\label{thm-codimension3k2}
     Let $\nu = a \sigma_3 + b \sigma_{2,1}$ be a nonzero, integral cohomology class in $G(2,n)$ with $n \geq 6$. Then $\nu$ is realizable over $\ZZ$ in $G(2, n)$ if and only if $a, b \geq 0$.
     \end{theorem}
\begin{proof}
  In $G(2,5)$, codimension 3 and dimension 3 classes agree. Hence, $a \sigma_3 + b \sigma_{2,1}$ is realizable over $\ZZ$ in $G(2,5)$ if $a,b >0$. By Propsition \ref{prop-bumpkn} (1), $I_{\ZZ}^3 (2,5) \subseteq I_{\ZZ}^3 (2,n)$ for $n >5$. Hence, $a \sigma_3 + b \sigma_{2,1}$ is realizable over $\ZZ$ in $G(2,n)$ if $a,b >0$. The classes $m \sigma_3$ and $m \sigma_{2,1}$ are realizable over $\ZZ$ for every $m > 0$ in $G(2,n)$ for $n>5$ by \cite[Theorem 1.1]{coskunrobles}. This concludes the proof. 
\end{proof}

\section{Dimension 4 realizable classes}

In this section, we discuss dimension 4 classes that are realizable.  For $G(3,6)$ and $G(2,n)$ we give a characterization of dimension 4 realizable classes, but emphasize that we do not obtain such a characterization of dimension 4 realizable classes for all Grassmannians.

\begin{proposition}
    Assume that $$a x^2y + b x^2 z + c xyz + d xz^2 + e yz^2$$ is realizable over $\ZZ$ in $\PP^2 \times \PP^1 \times \PP^2$. Then the following classes are realizable over $\ZZ$ in $G(4,8)$.
    \begin{enumerate}
        \item $a \sigma_{4,4,3} + b \sigma_{4,4,2,1} + c \sigma_{4,3,3,1} + d \sigma_{4,3,2,2} + e \sigma_{3,3,3,2}$.
        \item $a \sigma_{4,4,4} + (a+b+c) \sigma_{4,4,3,1} + (b+d) \sigma_{4,4,2,2} + (c+d+e) \sigma_{4,3,3,2}  + e \sigma_{3,3,3,3}$
    \end{enumerate}
\end{proposition}

\begin{proof}
    If the class in (1) is realizable over $\ZZ$, then by intersecting with a hyperplane class and using Pieri's formula we see that the class in (2) is realizable over $\ZZ$.

    Fix three transverse subspaces $A,B,C$ of $V$ of dimensions $3,2$ and $3$, respectively. Consider the incidence correspondence
    $$U := \{ (\ell, p_1, p_2, W) \ | \ \ell \subset A, p_1 \subset B, p_2 \subset C, \ell \subset W, \dim(W \cap \overline{p_1A}) \geq 3, W \subset \overline{p_2 AB} \}$$
    $$\subset \PP(A)^* \times \PP(B) \times \PP(C) \times G(4,8),$$
    parameterizing four-tuples of a $2$-dimensional subspace $\ell$ of $A$ and one-dimensional subspaces $p_1$ and $p_2$ of $B$ and $C$, respectively, and a $4$-dimensional subspace $W$ that contains $\ell$ has at least a  three-dimensional intersection with the span $\overline{p_1A}$ and is contained in the span $\overline{p_2 A B}$. Under this correspondence using Proposition \ref{prop-inductivestatement} and Theorem \ref{thm-mainconstruction2ind}, we see that surface classes in $\PP^2 \times \PP^1 \times \PP^2$ map to the claimed classes in $G(4,8)$.
\end{proof}

\begin{theorem}\label{thm-4dimG36}
    Let $\nu = a \sigma^{3,1} + b \sigma^{2,2} + c \sigma^{2,1,1}$ be an integral, nonzero four-dimensional class in $G(3,6)$. Then $\nu$ is realizable over $\ZZ$ if and only if $a,b,c \geq 0$ and one of the following holds:
\begin{enumerate}
  \item $a>0$ or $c>0$; or
  \item $(a,b,c)=(0,1,0)$.
\end{enumerate}
\end{theorem}

\begin{proof}
    Let $axy+bxz+cyz$ be a nonzero, effective curve class in $\PP^1 \times \PP^1 \times \PP^1$, where $x,y,z$ denote the pullbacks of the hyperplane class from the three projections. If at least two of $a,b,c$ are positive, then we can consider three general maps $$f_i : \PP^1 \to \PP^1$$ of degrees $a,b$ and $c$, respectively. This induces a birational map $f: \PP^1 \to \PP^1 \times \PP^1 \times \PP^1$ whose image is a curve with class $axy+bxz+cyz$.

    Fix three  transverse two-dimensional linear space $A,B,C$ in the underlying vector space $V$. 
    Consider the incidence correspondence
    $$U := \{ (p,q,r, W) \ | \ p\subset A, q\subset B, r\subset C, p \subset W, \dim(W \cap \overline{qA}) \geq 2, W \subset \overline{rAB}\}$$
    $$\subset \PP^1 \times \PP^1 \times \PP^1 \times G(3,6).$$
By Theorem \ref{thm-mainconstruction2ind}, an irreducible curve with class $axy+bxz +cyz$ in $\PP^1 \times \PP^1 \times \PP^1$ gives rise to an irreducible fourfold in $G(3,6)$ with class $a \sigma^{3,1} + b \sigma^{2,2} + c \sigma^{2,1,1}$.
Hence, such a class can be represented by an irreducible fourfold if $a,b,c \geq 0$ and at most one of the coefficients is 0.

There remains to discuss the case when only one of $a,b$ or $c$ is positive. By \cite[Theorem 1.1]{coskunrobles}, $a \sigma^{3,1}$ and $c \sigma^{2,1,1}$ are realizable over $\ZZ$ for every positive integer $a$ and $c$. By \cite{hong:rigidity}, $b\sigma^{2,2}$ is realizable over $\ZZ$ if and only if $b=1$. This concludes the proof of the theorem. 
\end{proof}

 We now turn to $G(2,n)$. In $G(2,5)$ a four-dimensional class has codimension 2, hence has been discussed in Theorem \ref{thm-codimension2}. We may therefore assume that $n\geq 6$. 

\begin{theorem}\label{thm-dim4G(2,n)}
  Let $\nu = a \sigma^4 + b \sigma^{3,1} + c \sigma^{2,2}$ be an integral, nonzero, four-dimensional class in $G(2,n)$ for $n \geq 6$. Then $\nu$ is realizable over $\ZZ$ if and only if $a,b,c \geq 0$, $b^2 \geq ac$ and one of the following holds:
\begin{enumerate}
    \item If $n>6$, we have $b>0$, or $b=c=0$, or $(a,b,c)=(0,0,1)$.
    \item If $n=6$, we have $b>0$, or $(a,b,c)=(1,0,0)$, or $(a,b,c)=(0,0,1)$.
   \end{enumerate}  
  \end{theorem}
\begin{proof}
We begin by considering the case of $G(2,6)$.  Fix two transverse three-dimensional subspaces $A,B$ in the underlying vector space $W$. Consider the incidence correspondence
$$U := \{ (p,q,W) \ | \ p \subset A, q\subset B, p \subset W \subset \overline{qA} \} \subset \PP^2 \times \PP^2 \times G(2,6),$$ where $p$ and $q$ are one-dimensional subspaces of $A$ and $B$, respectively, and $W$ is a two-dimensional subspace containing $p$ and contained in the span of $q$ and $A$. Given an irreducible surface with class $a x^2 + bxy + cy^2$, by Theorem \ref{thm-mainconstruction2}, we obtain an irreducible fourfold in $G(2,6)$ with class $a \sigma^4 + b \sigma^{3,1} + c \sigma^{2,2}$.  By \cite[Theorem 1]{huh2}, a surface class $a x^2 + bxy + cy^2$ is realizable when $b>0$ and $b^2 \geq ac$.

Suppose $a \sigma_4 + b \sigma_{3,1} + c \sigma_{2,2}$ is the class of an irreducible fourfold $Y$ in $G(2,6)$. We can embed this fourfold $Y$ in $G(3,7)$ by linearly embedding $V_6 \subset V_7$ and taking a vector $e$ not in $V_6$ and considering $\overline{eW}$ for $W \in Y$. The class of $Y$ in $G(3,7)$ is then $a \sigma_{4,4} + b \sigma_{4,3,1} + c \sigma_{4,2,2}$. Consider the birational map $f : \PP^4 \times G(2,6) \dashrightarrow G(3,7)$. We claim that we have the following pullbacks:
\begin{align*}
&f^* \sigma_{4,4} \ \ = \sigma_{4,4} + x \sigma_{4,3} + x^2 \sigma_{4,2} + x^3 \sigma_{4,1} + x^4 \sigma_4 \\
& f^* \sigma_{4,3,1}= x \sigma_{4,3} + x^2 \sigma_{4,2} + x^2 \sigma_{3,3} + x^3 \sigma_{4,1} + x^3 \sigma_{3,2} + x^4 \sigma_{3,1} \\
& f^* \sigma_{4,2,2} = x^2 \sigma_{4,2} + x^3 \sigma_{3,2} + x^4 \sigma_{2,2} 
\end{align*}
These formulae can be easily verified by Schubert calculus. The dual of a class $x^i \sigma_{\lambda_1, \lambda_2}$ is the Schubert class $x^{4-i} \sigma_{4 - \lambda_2, 4 - \lambda_1}$. The coefficient of $f^* \sigma_{\mu}$ is then given by $\sigma_{\mu} \cdot \sigma_{4- 
\lambda_2, 4- \lambda_1} \cdot \sigma_{4-i}.$ This product is easily computed by Pieri's rule.

Hence, if we intersect the class of $Y$ with $x^2$, we obtain a surface class
$$a x^2 \sigma_{4,4} + (a+b) x^3 \sigma_{4,3} + (a+b+c) x^4 \sigma_{4,2} + b x^4 \sigma_{3,3}. $$
The Hodge index inequality with respect to $x$ and $\sigma_1$ gives us
$$a(a+2b+c) \leq (a+b)^2,$$
which after simplifying becomes
$$ac\leq b^2.$$
Hence, we conclude that $b^2 \geq ac$ is a necessary inequality for the class to be realizable over $\ZZ$ or $\QQ$.

When $n=6$, the classes $a\sigma^4$ (resp. $c \sigma^{2,2}$) are realizable over $\ZZ$ if and only if $a=1$ (resp. $c=1$) by \cite{hong:rigidity}. This concludes the proof of (2).

    If $n > 6$, by Proposition \ref{prop-bumpkn} (1), we have $I_4^{\ZZ} (2,6) \subseteq I_4^{\ZZ} (2,n)$. Hence, any class with $b>0$ and $b^2 \geq ac$ is realizable over $\ZZ$. On the other hand, since $I_4^{\QQ}(2,n) \subseteq I_4^{\QQ}(2,6)$ by repeated applications of Proposition \ref{prop-dumpdown} (3), any realizable class must satisfy $b^2 \geq ac$. Hence, the only cases left to consider are when $b=0$. In that case, since $0 \geq ac \geq 0$, we must have either $a=0$ or $c=0$. By \cite[Theorem 1.1]{coskunrobles}, $a\sigma^4$ is realizable over $\ZZ$ for every positive $a$. By \cite{hong:rigidity}, $c \sigma^{2,2}$ is realizable over $\ZZ$ if and only if $c=1$. This concludes the prooof of (1).
\end{proof}

\section{Miscellaneous Results}

\subsection{Realizability in $G(2,n)$} In this subsection, we classify all cohomology classes in $G(2,n)$ that are realizable over $\QQ$.

\begin{theorem}\label{thm-g2ningeneral}
    Assume that $m < 2(n-2)$. Then a nonzero cohomology class $$\nu = \sum_{i=\max(0,m-n+2)}^{\left\lfloor \frac{m}{2} \right\rfloor} a_i \sigma^{m-i, i}$$ in $G(2,n)$ is realizable over $\QQ$ if and only if the coefficients $a_i$ form a log concave sequence of nonnegative rational numbers with no internal zeros.
\end{theorem}
\begin{proof}
Set $$m_1 = \left\lceil \frac{m}{2} \right\rceil \quad \mbox{and} \quad m_2 = \left\lfloor \frac{m}{2} \right\rfloor.$$
By Proposition \ref{prop-bundleoversegre}, the Schubert variety $\Sigma_{n-2-m_2} \subset G(2,n)$ is birational to a $\PP^{m_2}$-bundle over $\PP^{m_2} \times \PP^{n-m_2}$. By Theorem \ref{thm-mainconstruction2}, given  an irreducible  $m_1$-dimensional subvariety of $\PP^{m_2} \times \PP^{n-2-m_2}$ with class $\sum_{i= \max(0, m+2-n)}^{\min(m_2, n-2-m_2)} a_i x^{m_1-i}y^{n-m-2+i}$, there is an irreducible $m$-dimensional subvariety of $G(2,n)$ with class $\sum a_i \sigma^{m-i,i}$. By \cite[Theorem 21]{huh}, the class  $\sum a_i x^{m_1-i}y^{n-m-2+i}$ is realizable over $\QQ$ if and only if the coefficients $a_i$ form a nonzero log-concave sequence of nonnegative integers with no internal zeros. We conclude that the classes listed in the theorem are realizable over $\QQ$ in $G(2,n)$.

There remains to show that if the class $\nu$ is realizable over $\QQ$, then the coefficients must form a log concave sequence of nonnegative integers with no internal zeros. We will do this by induction on $m$. By Theorems \ref{thm-dimension2}, \ref{thm-codimension3k2} and \ref{thm-dim4G(2,n)}, this holds when $m\leq 4$. By Corollary \ref{cor-equalitydim} (1), we may assume that $n \leq m+2$. Given an irreducible subvariety $Y \subset G(2,n)$ of dimension $m$ with class $\nu = \sum_{i=\max(0,m-n+2)}^{m_2} a_i \sigma^{m-i, i}$, we can intersect it with a Schubert variety $\Sigma_{1,1}$ to get an irreducible subvariety of dimension $m-2$ in $G(2, n-1)$. 
By elementary Schubert calculus, the class of this variety is $\nu' = \sum_{i=\max(1,m-n+3)}^{m_2} a_i \sigma^{m-i-1, i-1}$.  By induction on $m$, if $i > \max(0, m-n+2)$, the coefficients $a_i$ form a log concave sequence of nonnegative integers with no internal zeros. There remains to show that the first three coefficients are also log concave. 

Embed the irreducible variety $Y$ with class $\nu$ in $G(3, n+1)$ to obtain an irreducible variety, which we will still call $Y$, with class $\sum_{i=\max(0,m-n+2)}^{m_2} a_i \sigma^{m-i, i}$. Consider the birational map $$f: \PP^{n-2} \times G(2, n) \dashrightarrow G(3,n+1)$$ obtained by fixing an $(n-1)$-dimensional subspace $\Lambda$ and a transverse $n$-dimensional subspace $\Omega$ and sending a $1$-dimensional subspace $W' \subset \Lambda$ and a two-dimensional subspace $W'' \subset \Omega$ to their span $\overline{W'W''}$. We consider $x^{n-4}\cdot [f^{-1}(Y)]$ and use the Hodge index theorem. We have the following relations:
\begin{align*}
   & x^{n-4} \cdot f^* \sigma^{n-2, m-n+2} = x^{n-2} \sigma^{2} + x^{n-3} \sigma^1 + x^{n-4}\sigma^0 \\
   & x^{n-4} \cdot f^* \sigma^{n-3, m-n+3} = x^{n-2} \sigma^2 + x^{n-2} \sigma^{1.1} + x^{n-3} \sigma^1 \\
   & x^{n-4} \cdot f^* \sigma^{n-4, m-n+4} = x^{n-2} \sigma^{1,1} \\
   & x^{n-4} \cdot f^* \sigma^{n-i, m-n+i} = 0 \ \ \mbox{for} \ i>4.
\end{align*}
Hence, 
 $$x^{n-4}\cdot [f^{-1}(Y)]= a x^{n-4} \sigma^0 + (a+b) x^{n-3} \sigma^1 + (a+b) x^{n-2} \sigma^2 + (b+c) x^{n-2} \sigma^{1,1}.$$ 
 The intersection matrix with respect to $x$ and $\sigma_1$ takes the form
 $$\left( \begin{array}{cc}
 a & a+b \\
 a+b & a+2b+c
 \end{array}
 \right).$$
 Hence, we obtain the inequality
 $$a(a+2b+c) - (a+b)^2 = ac - b^2 \leq 0.$$
 This concludes the inductive step and the proof of the theorem.
\end{proof}

\subsection{Realizability in $G(3,6)$}\label{sec-g36}
In this subsection, we describe all the cohomology classes in $G(3,6)$ that are realizable over $\ZZ$. All effective curve and divisor classes are realizable over $\ZZ$. The dimension and codimension 2 realizable classes have been classified in Theorem \ref{thm-dimension2} and Theorem \ref{thm-codimension2}, respectively. The dimension 3 and codimension 3 realizable classes have been classified in Theorem \ref{thm-dimension3} and Theorem \ref{thm-codimension3}, respectively. Finally, the dimension 4 classes in $G(3,6)$ that are realizable over $\ZZ$ have been classified in Theorem \ref{thm-4dimG36}. This only leaves classifying dimension 5 classes in $G(3,6)$ that are realizable over $\ZZ$, which is accomplished in the next theorem.

\begin{theorem}\label{thm-5dimG36}
    Let $\nu = a \sigma_{3,1} + b \sigma_{2,2} + c \sigma_{2,1,1}$ be an integral, nonzero five-dimensional class in $G(3,6)$. Then $\nu$ is realizable over $\ZZ$ if and only if $a,b,c \geq 0$ and one of the following holds:
\begin{enumerate}
  \item $a>0$ or $c>0$; or
  \item $(a,b,c)=(0,1,0)$.
\end{enumerate}
\end{theorem}
\begin{proof}
        Let $ax+by+cz$ be a nonzero, effective surface class in $\PP^1 \times \PP^1 \times \PP^1$, where $x,y,z$ denote the pullbacks of the hyperplane class from the three projections. If $a,b,c>0$, then the corresponding divisor class is very ample, hence is realizable over $\ZZ$ by Bertini's Theorem. If exactly one of $a,b$ or $c$ is zero (say $c=0$), then there exists irreducible curves of class $ax+by$ in $\PP^1 \times \PP^1$. The inverse image of such a curve under the projection is a surface with class $ax+by$ in $\PP^1 \times \PP^1$.

    Fix three  transverse two-dimensional linear space $A,B,C$ in the underlying vector space $V$. 
    Consider the incidence correspondence
    $$U := \{ (p,q,r, W) \ | \ p\subset A, q\subset B, r\subset C, p \subset W, \dim(W \cap \overline{qA} \geq 2, W \subset \overline{rAB}\}$$
    $$\subset \PP^1 \times \PP^1 \times \PP^1 \times G(3,6).$$
By  Theorem \ref{thm-mainconstruction2ind}, an irreducible surface with class $ax+by +cz$ in $\PP^1 \times \PP^1 \times \PP^1$ gives rise to an irreducible fivefold in $G(3,6)$ with class $a \sigma_{3,1} + b \sigma_{2,2} + c \sigma_{2,1,1}$.
Hence, such a class can be represented by an irreducible fivefold if $a,b,c \geq 0$ and at most one of the coefficients is 0.

There remains to discuss the case when only one of $a,b$ or $c$ is positive. By \cite[Theorem 1.1]{coskunrobles}, $a \sigma_{3,1}$ and $c \sigma_{2,1,1}$ are realizable over $\ZZ$ for every positive integer $a$ and $c$. By \cite{hong:rigidity}, $b\sigma_{3,1,1}$ is realizable over $\ZZ$ if and only if $b=1$. This concludes the proof of the theorem. 
\end{proof}

\subsection{Realizable classes in certain faces}

One can always use the construction in Proposition \ref{prop-inductivestatement} and Theorem \ref{thm-mainconstruction2ind} to produce irreducible subvarieties in $G(k,n)$. In general, this process produces irreducible subvarieties contained in a face of the effective cone. We illustrate the process with the following proposition, which is an immediate corollary of Theorem \ref{thm-mainconstruction2ind}.

\begin{proposition}\label{prop-face}
   Fix a strictly decreasing partition $n-k \geq \lambda_1 > \lambda_2 > \cdots > \lambda_k \geq 0$. Set $\lambda_0 = n-k$ Let $\nu= \sum a_{\mu} \sigma_{\mu}$ be the class of a $((k-1)(n-k) - \sum_{i=1}^{k-1} \lambda_i +m)$-dimensional cycle in $G(k,n)$ such that $\lambda_i \leq \mu_i \leq \lambda_{i-1}$ for $1 \leq i \leq k$. Then $\nu$ is realizable over $\ZZ$ if the $m$-dimensional class $\sum a_{\mu} \prod_{i=1}^k x_i^{\mu_i - \lambda_i}$ is realizable over $\ZZ$ in $\prod_{i=1}^k \PP^{\lambda_{i-1} - \lambda_i}$, where $x_i$ is the pullback of the hyperplane class by the $i$th projection.
\end{proposition}

\begin{remark}
By \cite[Theorem 1.8]{hhmww}, a surface class in a product of projective spaces is realizable over $\QQ$ if (and only if) the $(k\times k)$-symmetric intersection matrix with respect to the pullbacks of the hyperplane classes $x_i$ is Lorentzian. By Proposition \ref{prop-face}, setting $m=2$ we conclude that if the intersection matrix for the classes $x_i$ on the surface class $\sum a_{\mu} \prod_{i=1}^k x_i^{\mu_i - \lambda_i}$ in $\prod_{i=1}^k \PP^{\lambda_{i-1} - \lambda_i}$ is Lorentzian, then $\nu= \sum a_{\mu} \sigma_{\mu}$ is realizable over $\QQ$ in $G(k,n)$.
\end{remark}

\bibliographystyle{plain}

\end{document}